\newif\iflatexe
\latexefalse

\iflatexe
\documentclass{article}
  \usepackage{latexsym,a4,amscd}
\else
    \documentstyle[a4,amscd]{article}
\fi
\makeatletter
\@addtoreset{equation}{subsection}

\makeatother

\newcommand{\separate}{\medskip\noindent}

\def\newsection{ \separate
   \refstepcounter{subsection} 
   {\large\bf \thesubsection\kern.3em}
}
\def\mytheorem#1{
   \separate{\large\bf Theorem#1:\kern.3em}} 

\def\mylemma#1{
   \separate{\large\bf Lemma#1:\kern.3em}} 

\def\mycorollary#1{
   \separate{\large\bf Corollary#1:\kern.3em}} 

\def\myprop#1{
   \separate{\large\bf Proposition#1:\kern.3em}} 

\def\myremark#1{
   \separate{\large\bf Remark#1:\kern.3em}} 

\def\mydefinition#1{
   \separate{\large\bf Definition#1:\kern.3em}} 

\def\myconjecture#1{
   \separate{\large\bf Conjecture#1:\kern.3em}} 

\def\Proof{\separate\underline{Proof:}\kern1em}

\newcommand{\GANZ}{{\sf Z\hspace*{-0.4em}Z}}
\newcommand{\REELL}{{\setlength{\unitlength}{1em}
                     \begin{picture}(0.75,1)
                     \put(0,0){\line(0,1){0.69}}
                     \put(0,0){\sf R}
                     \end{picture}
                   }}
\newcommand{\KOMPLEX}{{\setlength{\unitlength}{1em}
                     \begin{picture}(0.7,1)
                     \put(0.34,0){\line(0,1){0.65}}
                     \put(0,0){\sf C}
                     \end{picture}
                   }}
\newcommand{\NATUR}{{\setlength{\unitlength}{1em}
                     \begin{picture}(0.75,1)
                     \put(0,0){\line(0,1){0.69}}
                     \put(0,0){\sf N}
                     \end{picture}
                   }}
\newcommand{\FGANZ}{\mbox{\tiny{\rm Z\hspace*{-0.45em}Z}}}
\newcommand{\FNATUR}{\mbox{\tiny{\setlength{\unitlength}{1em}
                     \begin{picture}(0.6,0.5)
                     \put(0,0){\line(0,1){0.48}}
                     \put(0,0){\rm N}
                     \end{picture}
                   }}}
\newcommand{\FKOMPLEX}{\mbox{\tiny{\setlength{\unitlength}{1em}
                               \begin{picture}(0.6,0.5)
                               \put(0.34,0){\line(0,1){0.47}}
                               \put(0,0){\rm C}
                               \end{picture}
                              }}}
\newcommand{\FREELL}{\mbox{\tiny{\setlength{\unitlength}{1em}
                     \begin{picture}(0.6,0.5)
                     \put(0,0){\line(0,1){0.48}}
                     \put(0,0){\rm R}
                     \end{picture}
                   }}}
\newcommand{\Brr}{{\mathchoice{\REELL}{\REELL}{\!\FREELL}{\!\FREELL}}}
\newcommand{\Bcc}{{\mathchoice{\KOMPLEX}{\KOMPLEX}{\!\FKOMPLEX}{
\!\FKOMPLEX}}}
\newcommand{\Bnn}{{\mathchoice{\NATUR}{\NATUR}{\FNATUR}{\FNATUR}}}
\newcommand{\Bii}{{\mathchoice{\GANZ}{\GANZ}{\FGANZ}{\FGANZ}}}
\newcommand{\CPE}{{\mathchoice{\Bcc{\rm P}_1}{\Bcc{\rm P}_1}{\Bcc\!\!
\mbox{\rm\tiny P}_1}{\Bcc\!\!\mbox{\rm\tiny P}_1}}}
\newcommand{\ClPE}{{\mathchoice{\Bcc_\lambda{\rm P}_1}{
\Bcc_\lambda{\rm P}_1}{\Bcc_\lambda\!\!
\mbox{\rm\tiny P}_1}{\Bcc_\lambda\!\!\mbox{\rm\tiny P}_1}}}
\newcommand{\CnPE}{{\mathchoice{\Bcc_\nu{\rm P}_1}{\Bcc_\nu{\rm P}_1}{
\Bcc_\nu\!\!
\mbox{\rm\tiny P}_1}{\Bcc_\nu\!\!\mbox{\rm\tiny P}_1}}}
\newcommand{\ad}{{\mathchoice{\mbox{\rm ad}}{\mbox{\rm ad}}{%
\mbox{\scriptsize\rm ad}}{\mbox{\scriptsize\rm ad}}}}

\def\QED{\hfill$\Box$}
\def\inv{^{-1}}
\def\BEA{\begin{eqnarray}}
\def\EEA{\end{eqnarray}}
\def\BEQ{\begin{equation}}
\def\EEQ{\end{equation}}
\def\bref#1{(\ref{#1})}
\def\tmatrix#1#2#3#4{
    \left(\begin{array}{cc} #1 & #2 \\ #3 & #4 \end{array}\right)}
\def\LieSL{{\mbox{\bf SL}}}
\def\LieSO{{\mbox{\bf SO}}}

\def\LieSU{{\mbox{\bf SU}}}
\def\LieB{{\mbox{\bf B}}}

\def\LieU{{\mbox{\bf U}}}
\def\Liesl{{\mbox{\bf sl}}}
\def\Liesu{{\mbox{\bf su}}}

\def\dprime{{\prime\prime}}

\def\vecr{{\mbox{\bf r}}}
\def\vecsigma{{\sigma\kern-2.25mm\sigma}}

\def\tb{\tilde{b}}

\def\tf{\tilde{f}}

\def\tw{\tilde{w}}
\def\tmu{\tilde{\mu}}

\def\tsigma{\tilde{\sigma}}
\def\FCinvolution{\hat{\sigma}}

\def\FA{{\cal A}}
\def\FC{{\cal C}}
\def\FCstar{{{\cal C}^\ast}}
\def\FD{{\cal D}}

\def\FG{{\cal G}}

\def\FM{{\cal M}}

\def\FO{{\cal O}}

\def\FRL{{{\cal R}\kern1mm_\lambda}}
\def\FS{{\cal S}}
\def\FT{{\cal T}}
\def\FU{{\cal U}}
\def\FZ{{\cal Z}}

\def\der#1#2{{{\rm d}#1\over {\rm d}#2}}

\def\diff{{\rm d}}

\def\diag{{\mbox{\rm diag}}}
\def\Ker{{\mbox{\rm Ker}}}

\def\Aut{{\mbox{\rm Aut}}}
\def\Sym{{\mbox{\rm Sym}}}
\def\Per{{\mbox{\rm Per}}}
\def\OAff{{\mbox{\rm OAff}}}

\def\twospace{{\Brr^{\!\lower2pt\hbox{\mbox{\rm\scriptsize{2}}}}}}
\def\threespace{{\Brr^{\!\lower2pt\hbox{\mbox{\rm\scriptsize{3}}}}}}
\def\cstar{{\Bcc^{\!\lower2pt\hbox{\mbox{\scriptsize{$\ast$}}}}}}
\def\cnstar{{\Bcc_\nu^{\!\lower2pt\hbox{\mbox{\scriptsize{$\ast$}}}}}}
\def\clstar{{\Bcc_\lambda^{\!\lower2pt\hbox{\mbox{\scriptsize{$\ast$}}}}}}
\def\CPEtwo{{\CPE^{\!\!\lower2pt\hbox{\mbox{\rm\scriptsize{2}}}}}}
\def\Real{{\mbox{\rm Re}}}
\def\Imag{{\mbox{\rm Im}}}

\def\tr{{\mbox{\rm tr}}}

\def\hphi{\hat{\phi}}

\def\ha{\hat{a}}
\def\hb{\hat{b}}
\def\hc{\hat{c}}
\def\hd{\hat{d}}

\def\hS{\hat{S}}

\def\hrho{\hat{\rho}}
\def\ca{\check{a}}

\def\tp{\tilde{p}}

\def\tphi{\tilde{\phi}}
\def\Equer{{\overline{E}}}

\def\qquer{{\overline{q}}}
\def\zquer{{\overline{z}}}

\def\lambdaquer{\overline{\lambda}}
\def\nuquer{\overline{\nu}}
\def\pder#1#2{{\partial #1\over\partial #2}}

\begin{document}

\renewcommand{\thefootnote}{\fnsymbol{footnote}}
\begin{center}
{\LARGE On constant mean curvature surfaces with periodic metric}

\vskip1cm
\begin{minipage}{6cm}
\begin{center}
J.~Dorfmeister\\
Department of Mathematics\\ University of Kansas\\ Lawrence, KS 66045
\end{center}
\end{minipage}
\begin{minipage}{6cm}
\begin{center}
G.~Haak\footnotemark[1]\\ Fachbereich Mathematik\\ TU Berlin\\ D-10623 Berlin
\end{center}
\end{minipage}
\vspace{0.5cm}

\end{center}

\footnotetext[1]{supported by Sonderforschungsbereich 288}

\section{Introduction}\refstepcounter{subsection} \label{introduction}
\message{[Introduction]}
In recent years two independent approaches to the construction of 
CMC-immersions in $\threespace$ were developed. 
One of them, using finite type solutions of the $\sinh$-Gordon
equation, leads to the classification of CMC-tori
by Pinkall and Sterling~\cite{PiSt:1}.

While this approach allowed
A.~Bobenko~\cite{Bo:1} to describe CMC-tori explicitly in terms of 
theta functions, it turned out not to be easily adaptable to the more
general situation of CMC-surfaces with umbilics. 

The class of CMC-surfaces with umbilics contains 
the CMC-surfaces of higher genus, whose existence was proved by
N.~Kapouleas for arbitrary genus~\cite{Ka:1,Ka:2}. Therefore, it
proved to be necessary to find a more general description of CMC-immersions.

Another approach to CMC-surfaces was developed by F.~Pedit, H.~Wu and
one of the authors (J.D.). It is commonly called the DPW method. It gives a
description of all immersed CMC-surfaces in $\threespace$ with or without
umbilics. More specifically, it allows the construction of 
CMC-immersions from a meromorphic and a holomorphic function, for which
reason it can be called a Weierstra\ss\ type representation of 
CMC-immersions. These functions, one of which is the coefficient of
the Hopf differential of
the surface to be constructed, are only subject to some simple
restrictions \cite{DoHa:1}.

The main problem with the DPW method as with the Weierstra\ss\
representation for minimal surfaces is the fact, that in spite of the
simplicity of the initial data, it is not easy to extract
further geometric information out of these data. For one point, it turns
out to be very hard to impose restrictions on the symmetries of the
surface in $\threespace$ or on the Fuchsian or elementary group of the
underlying Riemann surface. One of the cases where this would be
necessary is the case of compact surfaces.
In fact, it was not even clear how to reproduce the results of Pinkall,
Sterling~\cite{PiSt:1} and Bobenko~\cite{Bo:1} with the DPW method.

In this paper we will solve this problem and reproduce the
classification of CMC-tori using extended frames without referring to
the metric, i.e., the $\sinh$-Gordon equation, explicitly.
This result will be proved in a more general framework, which was
developed in \cite{DoHa:2} by the authors for the investigation of 
CMC-surfaces with symmetries in $\threespace$. Here we will actually derive a
classification of all finite type surfaces with periodic metric, i.e., 
surfaces whose metric is a finite type solution of the $\sinh$-Gordon
equation which is invariant under a group of translations in
$\Bcc$. This classification is of course also implicit in \cite{PiSt:1}.

But first, let us define more precisely, what we mean by a CMC-surface
with periodic metric. Here we follow \cite{DoHa:2}.

\newsection \label{intro11}
Let $\Psi:\Bcc\rightarrow\threespace$ be a conformal CMC-immersion
of the complex plane into $\threespace$.
We further set $\OAff(\threespace)$ to be the group of proper (i.e.,
orientation preserving) Euclidean motions in $\threespace$.

\mydefinition{}
Let $\Psi:\Bcc\rightarrow\threespace$ be defined as above. Then we define
the following subset $\Sym(\Psi)$ of $\Bcc$:
A complex constant $q$ is in $\Sym(\Psi)$, iff
there exists a proper Euclidean motion $T_q\in\OAff(\threespace)$, s.t.\
\BEQ \label{symmetrydef}
\Psi(z+q)=T_q(\Psi(z))\kern1cm\mbox{\rm for all $z\in\Bcc$}.
\EEQ
If $\Sym(\Psi)\neq\{0\}$, then we say that $\Psi$ is a surface with 
periodic metric.

We further define the subset $\Per(\Psi)$ of $\Sym(\Psi)$ by
\BEQ \label{perdef}
\Per(\Psi)=
\{q\in\Bcc| \Psi(z+q)=\Psi(z),\;\mbox{\rm for all $z\in\Bcc$}\}.
\EEQ
If $\Per(\Psi)\neq\{0\}$, then we say that $\Psi$ is a periodic surface.

\myremark{} 1.
For the investigation of translational symmetries, it is of course enough to
consider only those Riemannian surfaces whose universal cover is the whole
complex plane. Therefore, we will not introduce the more complex
notation of \cite{DoHa:2}. With the exception of the following theorem
we will always set 
$M=\FD=\Bcc$ and $\Phi=\Psi$ (see \cite[Sec.~2]{DoHa:2}).
From time to time we will also identify the elements of $\Sym(\Psi)$
and $\Per(\Psi)$ with the corresponding translations, i.e., we will
identify $\Sym(\Psi)$ and $\Per(\Psi)$ with subgroups of $\Aut\Bcc$,
the set of biholomorphic automorphisms of $\Bcc$.

2. To further compare our notation with that of \cite{DoHa:2}, we
note, that $\Sym(\Psi)$ can be identified isomorphically with the
subgroup of translations in $\Aut_\Psi\Bcc$. The group of all
automorphisms in $\Aut\Bcc$, for which the immersion $\Psi$
satisfies a relation of type~\bref{symmetrydef}, is defined in
\cite[Def.~2.6]{DoHa:2}. It is a closed Lie-subgroup of $\Aut\Bcc$. 
If we also introduce the homomorphism
$\psi:\Aut_\Psi\Bcc\rightarrow\OAff(\threespace)$ as in
\cite[Def.~2.9]{DoHa:2}, then $\psi$ is analytic
\cite[Theorem~2.9]{DoHa:2} and its image $\Aut\Psi(\Bcc)$ is, by
\cite[Def.~2.6 and Cor.~2.10]{DoHa:2}, the set of all proper Euclidean motions
which preserve the image $\Psi(\Bcc)$. Moreover, $Ker\psi\supset\Per(\Psi)$. 
Since by \cite[Lemma~2.9]{DoHa:2} $\Ker\psi$ acts freely on $\Bcc$, we have
that $\Ker\psi$ consists only of translations, whence $\Ker\psi=\Per(\Psi)$.

\separate
We will collect some facts about the groups $\Sym(\Psi)$ and
$\Per(\Psi)$ which follow from results in \cite{DoHa:2}.
For the definition of the associated family in this context see
\cite[Sect.~2.5]{DoHa:2}.

\mytheorem{} {\em
Let $(\Bcc,\Psi)$ be a conformal CMC-immersion. Then the following holds:

1. If $\Sym(\Psi)$ is nondiscrete, then $\Psi(\Bcc)$ is in the associated
family of a Delaunay surface.

2. If $\Sym(\Psi)$ is discrete and has a nonzero element, 
then its elements form a lattice in $\Bcc$ with one or two generators.

3. $\Per(\Psi)$ is always discrete and, if nontrivial, forms a lattice
in $\Bcc$ with one or two generators.

4. If $\Per(\Psi)\neq\{0\}$, then there exists a
Riemannian surface $M$ with universal cover $\pi:\Bcc\rightarrow M$ 
and a CMC-immersion $\Phi:M\rightarrow\threespace$,
s.t.\ $\Gamma=\Per(\Psi)$ is the elementary 
group of $M$ and $\Psi=\Phi\circ\pi$.

5. In particular, if $\Per(\Psi)$ forms a lattice with 
two generators, then $\Psi(\Bcc)$ is a CMC-torus.
}

\Proof
We use Remark~\ref{intro11} to compare with \cite{DoHa:2}.
If $\Sym(\Psi)$ is nondiscrete, then so is $\Aut_\Psi(\Bcc)$. Since
$\Aut_\Psi\Bcc$ is closed, it contains a one-parameter subgroup, which
acts, by~\cite[Cor.~2.6]{DoHa:2}, as a one-parameter group of 
self-isometries for $\Psi$. We
can therefore use \cite[Theorem~2.14 and Lemma~2.15]{DoHa:2} to prove
the first statement of the theorem.

From the definition it follows immediately, that with $q\in\Sym(\Psi)$ also
$nq\in\Sym(\Psi)$ for all $n\in\Bii$ and $T_{nq}=T_q^n$.
This implies, that the elements of $\Sym(\Psi)$ form a lattice,
which can have no more than two generators in $\Bcc\cong\twospace$.
This proves the second statement.

The third and fourth statement follows with Remark~\ref{intro11}
from \cite[Proposition~2.11]{DoHa:2}. For the definition of a
CMC-immersion of a Riemann manifold into $\threespace$
see~\cite[Sect.~2.1,2.2]{DoHa:2}.

Finally, if $\Per(\Psi)$ is a two-dimensional lattice, then
$(M,\Phi)$ is a CMC-torus, which finishes the proof.
\QED

\separate
The paper starts in Section~\ref{DPW} with a short description of the
DPW method for later comparison with other formulations. We will also
give the transformation properties of an extended frame under a
translation in $\Sym(\Psi)$ as they were derived in
\cite[Sec.~3.4]{DoHa:2}. In this section we will work with loop groups
on arbitrary circles $C_r$, $0<r<1$, which allows for a more general
definition of the dressing action (see also~\cite{BuPe:1,DoMcPeWu:1}).

In Chapter~\ref{transsym} we consider surfaces obtained by
$r$-dressing~\cite{DoWu:1,BuPe:1} from the standard cylinder:
$$
h_+(\lambda)e^{(\lambda\inv z-\lambda\zquer)A}
=F(z,\zquer,\lambda)p_+(z,\zquer,\lambda).
$$
Here $A=\tmatrix0110$ and $F$ is the (unitary) frame of the surface,
while $p_+$ and $h_+$, in a Fourier expansion, does not ontain any
negative powers of $\lambda$. By a result of~\cite{DoWu:1} (see
also~\cite{BuPe:1}) all CMC-tori can be obtained this way. The main
result of this chapter considers translations $z\mapsto z+q$ and
gives necessary and sufficient conditions on $h_+$ such that the
frame $F$ at $z+q$ is a $z$-independent rotation of $F$ at $z$ for all
$z\in\Bcc$ (Theorem~\ref{nectheorem}, Theorem~\ref{sufftheorem}). This
is equivalent with the metric of the surface being periodic with
period $q$. As a consequence, in Theorem~\ref{spaceclosing} we give
necessary and sufficient conditions, that a CMC-surface, obtained by
$r$-dressing of the standard cylinder by $h_+$, admits the translation
$z\mapsto z+q$ in its fundamental group. This describes all
topological cylinders obtained by $r$-dressing the standard cylinder.

The results of Chapter~3 are phrased in terms of certain even rational
functions $a^2(\lambda^2),b^2(\lambda^2),c^2(\lambda^2)$ and
meromorphic functions $\alpha,\beta^2$ on $\Bcc\setminus\{0\}$. 
In Chapter~\ref{periodicsurfaces} we consider the hyperelliptic
surface $\FC$ with branchpoints exactly at the zeroes and poles of odd
order of $a^2$ and a natural two-fold covering $\FC^\prime$ of
$\FC$. It turns out that all functions involved in the $r$-dressing,
and a priori only defined on the disk of radius $r$ around
$\lambda=0$, can be extended holomorphically to the
hyperelliptic surface $\FC$ or at least to the hyperelliptic surface
$\FC^\prime$. As an application of this extension we show
(Theorem~\ref{finitetype}), that {\em every CMC-surface with periodic
metric in the $r$-dressing orbit of the standard cylinder is of finite
type} (for a definition of ``finite type'' in this framework see
\cite{BuPe:1,DoWu:1,PiSt:1}).

In Chapter~\ref{tori} we start from a hyperelliptic surface $\FC$ and
define $\alpha,\beta^2,a^2,b^2,c^2$ in a natural way, thus producing
$a^2$ with only simple poles. In this case, which is generic for
CMC-tori, we reproduce the dressing
matrix and thus are able to describe (generically) all CMC-tori 
in the $r$-dressing orbit of the cylinder
(Section~\ref{bcconstruction}).
If $g$ is the genus of $\FC$, then our construction yields for every
CMC-torus in the $r$-dressing orbit of the cylinder naturally a
$g-1$-parameter family of tori. This fits well to the purely
algebro-geometric construction of such families of tori as given
in~\cite{Bo:1,Ja:1}.

It would be interesting to expand our construction in
Chapter~\ref{tori} to the ``non-generic'' case, i.e., to the case
where $a^2$ does not only have simple poles and zeroes of odd order.

\section{The DPW method} \label{DPW} \message{[DPW]}
Let us recall for the reader's convenience the cornerstones of the DPW method.
For a more detailed reference the reader should consult
\cite{DoPeWu:1,DoMcPeWu:1} and the appendix of \cite{DoHa:1}.
In this section we will follow
the conventions of \cite{DoHa:1}. By a CMC-immersion we will always
understand a conformal immersion of constant mean curvature. For a
justification of this restriction, see~\cite[Section~2]{DoHa:2}.

\newsection{} \label{DPW1}
Let $\Psi:\FD\rightarrow\threespace$ be a CMC-immersion.
Here, $\FD\subset\Bcc$ is the open unit disk or the
whole complex plane. We take $\Psi$ as a conformal chart and the
metric to be $\diff s^2=e^u(\diff x^2+\diff y^2)$, $u:\FD\rightarrow\Brr$. 
Then
\BEQ
\FU=(e^{-{u\over2}}\Psi_x,e^{-{u\over2}}\Psi_y,N):\FD\longrightarrow\LieSO(3)
\EEQ
is an orthonormal frame which we normalize by
\BEQ \label{FUinitial}
\FU(x=0,y=0)=I.
\EEQ

In complex coordinates $z=x+iy$, $\zquer=x-iy$, we have
\BEQ
\langle\Psi_z,\Psi_z\rangle=\langle\Psi_\zquer,\Psi_\zquer\rangle=0,\kern3cm
\langle\Psi_z,\Psi_\zquer\rangle=\frac{1}{2}e^u.
\EEQ
Using the definitions
\BEQ
E=\langle\Psi_{zz},N\rangle,\kern2cm H=2e^{-u}\langle\Psi_{z\zquer},N\rangle
\EEQ
we get for the second fundamental form
\BEQ
II=\frac{1}{2}\tmatrix{(E+\Equer)+He^u}{i(E-\Equer)}{i(E-\Equer)}{
-(E+\Equer)+He^u}.
\EEQ
Thus, $H$ is the mean curvature and therefore constant.
The Gau\ss-Codazzi equations take the form
\BEA \label{GaussCodazzi1}
u_{z\zquer}+{1\over2}e^uH^2-2e^{-u}|E|^2 & = & 0, \\
E_\zquer={1\over2}e^uH_z=0,\label{GaussCodazzi2}
\EEA
from the second of which it follows, that $E\diff z^2$ is a holomorphic
differential, the Hopf differential.
For $|E|=1$ and $H=2$ the first equation becomes the $\sinh$-Gordon equation.

By using the spinor representation
$J:\vecr\mapsto-\frac{i}{2}\vecr\vecsigma$, where $\vecsigma$ is the
vector, whose components are the Pauli matrices, we identify vectors in
$\threespace$ with matrices in $\Liesu(2)$. This induces an identification
of orthonormal frames with matrices in $\LieSU(2)$ (see \cite{DoHa:1})
which is unique up to the multiplication with $-I$ in
$\LieSU(2)$.
If we identify the frame $\FU$ in this way with a map
\BEQ
P:\FD\rightarrow\LieSU(2),
\EEQ
then the initial condition~\bref{FUinitial} is compatible with
\BEQ \label{Pinitial}
P(0)=I,
\EEQ
which fixes the lift $P$ uniquely.
If we define
\BEA
U=P\inv P_z=\tmatrix{-{1\over4}u_z}{Ee^{-{u\over2}}}{
-{1\over2}e^{u\over2}H}{{1\over4}u_z},\\
V=P\inv P_\zquer=\tmatrix{{1\over4}u_\zquer}{{1\over2}He^{u\over2}}{
-\Equer e^{-{u\over2}}}{-{1\over4}u_\zquer},
\EEA
then the integrability condition
\BEQ \label{ZCC}
U_\zquer-V_z-[U,V]=0
\EEQ
is equivalent to the Gau\ss-Codazzi equations \bref{GaussCodazzi1}
and~\bref{GaussCodazzi2}.

By the substitution $E\rightarrow \lambda^{-2}E$, $\lambda\in S^1$,
which doesn't change the Gau\ss-Codazzi equations, we introduce the
so called spectral parameter $\lambda$. After a simple
re-gauging of the frames (for details see the appendix of~\cite{DoHa:1}) 
we get the following form of the now
$\lambda$-dependent matrices $U(\lambda)$ and $V(\lambda)$:
\BEA
U(\lambda) & = & \tmatrix{{1\over4}u_z}{
-{1\over2}\lambda\inv He^{u\over2}}{
\lambda\inv Ee^{-{u\over2}}}{-{1\over4}u_z}, \label{UDEF} \\
V(\lambda) & = &
\tmatrix{-{1\over4}u_\zquer}{-\lambda\Equer e^{-{u\over2}}}{
{1\over2}\lambda He^{u\over2}}{{1\over4}u_\zquer}.\label{VDEF}
\EEA
By integration, i.e., solving
\BEA
F\inv F_z & = & U,\\
F\inv F_\zquer & = & V, \label{VF}
\EEA
using the initial condition
\BEQ \label{Finitial}
F(0,\lambda)=I,\kern1cm\mbox{\rm for all $\lambda\in S^1$},
\EEQ
we get a $\lambda$-dependent frame
$F(z,\zquer,\lambda)$, which coincides at $\lambda=1$ with the old
frame $P(z,\zquer)$. We call this map from $\FD\times S^1$ to
$\LieSU(2)$ the {\em extended frame} of the immersion $\Psi$.
It should be noted, that by Eq.~\bref{UDEF} and~\bref{VDEF} the
differential of $F$ has the form
\BEQ \label{alphaform}
\alpha(z,\zquer,\lambda)=\lambda\inv\alpha_{-1}(z,\zquer)\diff z
+\alpha_0(z,\zquer)+\lambda\alpha_1(z,\zquer)\diff\zquer,\kern1cm
\alpha_1(z,\zquer)=-\overline{\alpha_{-1}(z,\zquer)}^\top.
\EEQ
This map is now interpreted as taking values in a certain loop
algebra $\Lambda_r\Liesu(2)_\sigma$.

\newsection{} \label{DPWloopgroups}
For each real constant $r$, $0<r<1$, 
let $\Lambda_r\LieSL(2,\Bcc)_\sigma$ denote the group of 
smooth maps $g(\lambda)$ from $C_r$, the circle of radius $r$, to
$\LieSL(2,\Bcc)$, which satisfy the twisting condition
\BEQ \label{DPWtwistcond}
g(-\lambda)=\sigma(g(\lambda)),
\EEQ
where $\sigma:\LieSL(2,\Bcc)\rightarrow\LieSL(2,\Bcc)$ 
is the group automorphism
of order $2$, which is given by conjugation with the Pauli matrix 
\BEQ
\sigma_3=\tmatrix100{-1}.
\EEQ
The Lie algebras of these groups, which we denote by
$\Lambda_r\Liesl(2,\Bcc)_\sigma$, consist of maps 
$x:C_r\rightarrow\Liesl(2,\Bcc)$, which satisfy a similar
twisting condition as the group elements
\BEQ
x(-\lambda)=\sigma_3 x(\lambda)\sigma_3.
\EEQ
In order to make these loop groups complex Banach Lie groups, we equip them,
as in \cite{DoPeWu:1}, with some $H^s$-topology for $s>{1\over2}$.
Elements of these twisted loop
groups are matrices with off-diagonal entries which are odd functions, and
diagonal entries which are even functions in the parameter $\lambda$.
All entries are in the Banach algebra $\FA_r$ of $H^s$-smooth
functions on $C_r$.

Furthermore, we will use the following subgroups of
$\Lambda_r\LieSL(2,\Bcc)_\sigma$: 
Let $\LieB$ be a subgroup of $\LieSL(2,\Bcc)$ and 
$\Lambda_{r,B}^+\LieSL(2,\Bcc)_\sigma$ be the group of maps in
$\Lambda_r\LieSL(2,\Bcc)_\sigma$, which can be extended to holomorphic maps on
\BEQ
I^{(r)}=\{\lambda\in\Bcc; |\lambda|<r\},
\EEQ
the interior of the circle $C_r$, and take values in $\LieB$ at $\lambda=0$.
Analogously, let $\Lambda_{r,B}^-\LieSL(2,\Bcc)_\sigma$ 
be the group of maps in $\Lambda_r\LieSL(2,\Bcc)_\sigma$, which can be extended
to the exterior 
\BEQ
E^{(r)}=\{\lambda\in\CPE;|\lambda|>r\}
\EEQ
of $C_r$ and take values in $\LieB$ at $\lambda=\infty$. 
If $\LieB=\{I\}$ (based loops) 
we write the subscript $\ast$ instead of $\LieB$, if
$\LieB=\LieSL(2,\Bcc)$ we omit the subscript for $\Lambda$ entirely.

Also, by an abuse of notation, we will denote by 
$\Lambda_r\LieSU(2)_\sigma$ the subgroup of maps in
$\Lambda_r\LieSL(2,\Bcc)_\sigma$, which can be extended holomorphically
to the open annulus
\BEQ
A^{(r)}=\{\lambda\in\Bcc;r<|\lambda|<\frac{1}{r}\}
\EEQ
and take values in $\LieSU(2)$ on the unit circle.

Corresponding to these subgroups, we analogously define Lie subalgebras of
$\Lambda_r\Liesl(2,\Bcc)_\sigma$.

We quote the following results from~\cite{Mc:1} and \cite{DoPeWu:1}:

\separate (i) For each solvable subgroup $\LieB$ of 
$\LieSL(2,\Bcc)$, which satisfies $\LieSU(2)\cdot\LieB=\LieSL(2,\Bcc)$ and
$\LieSU(2)\cap\LieB=\{I\}$, multiplication 
$$
\Lambda_r\LieSU(2)_\sigma\times\Lambda_{r,B}^+\LieSL(2,\Bcc)_\sigma
\longrightarrow\Lambda_r\LieSL(2,\Bcc)_\sigma
$$
is a diffeomorphism onto.
The associated splitting
\BEQ \label{Iwasawa}
g=F g_+
\EEQ
of an element $g$ of $\Lambda_r\LieSL(2,\Bcc)_\sigma$, s.t.\
$F\in\Lambda_r\LieSU(2)_\sigma$ and
$g_+\in\Lambda^+_{r,B}\LieSL(2,\Bcc)_\sigma$
will be called Iwasawa decomposition. In the following, we will fix the
group $\LieB$ as the group of upper triangular $2\times2$-matrices 
with real positive entries on the diagonal.

\separate (ii) Multiplication 
\BEQ \label{Birkhoff}
\Lambda^-_{r,\ast}\LieSL(2,\Bcc)_\sigma\times\Lambda_r^+\LieSL(2,\Bcc)_\sigma
\longrightarrow\Lambda_r\LieSL(2,\Bcc)_\sigma
\EEQ
is a diffeomorphism onto the open and dense subset 
$\Lambda^-_{r,\ast}\LieSL(2,\Bcc)_\sigma\cdot\Lambda_r^+\LieSL(2,\Bcc)_\sigma$
of $\Lambda_r\LieSL(2,\Bcc)_\sigma$, called the ``big cell'' \cite{SeWi:1}.
The associated splitting
\BEQ
g=g_-g_+
\EEQ
of an element $g$ of the big cell, where
$g_-\in\Lambda^-_{r,\ast}\LieSL(2,\Bcc)_\sigma$ and
$g_+\in\Lambda_r^+\LieSL(2,\Bcc)_\sigma$, 
will be called Birkhoff factorization.

\separate 
Clearly, for arbitrary $0<r<1$, 
the matrices $U(\lambda)$ and $V(\lambda)$ defined above take
values in $\Lambda_r\Liesu(2)_\sigma$, and so does $\alpha$.

We also get

\mylemma{} {\em
The extended frame $F(z,\zquer,\lambda)$ of a CMC-immersion can be
extended holomorphically in $\lambda$ to
$\cstar=\Bcc\setminus\{0\}$.
}

\Proof
For $n\in\Bnn$ define
$K_n=\{\lambda\in\Bcc;\frac{1}{n}\leq|\lambda|\leq n\}$. Then
$\{K_n\}_{n\in\Bnn}$ is a sequence of compact subsets of $\cstar$,
which exhaust $\cstar$. For each $n\in\Bnn$,
$\alpha(z,\zquer,\lambda)$, defined
by~\bref{alphaform}, takes values in the Banach algebra $C(K_n)$ of bounded
functions on $K_n$, which are holomorphic in the interior of $K_n$.
Since also the initial condition~\bref{Finitial} for $F$ is in
$C(K_n)$, we get that for each $n\in\Bnn$, the solution
$F(z,\zquer,\lambda)$ of the differential equation $F\inv\diff
F=\alpha$ takes values in $C(K_n)$.
Therefore, $F$ can be continued holomorphically in the $\lambda$-plane to
the open set
\BEQ
\bigcup_{n\in\Bnn}K_n=\cstar.
\EEQ
\mbox{}\QED

Using the initial condition $F(0,\lambda)=I$ we thus see, 
that, for each radius $r$, the extended frame
$F$ is, by restriction to $C_r$, a map from $\FD$ to 
$\Lambda_r\LieSU(2)_\sigma$. We will use
this fact to identify the extended frames with elements of
$\Lambda_r\LieSU(2)_\sigma$.

\separate
If $F(z,\zquer,\lambda)$ is the extended frame of a CMC-immersion
$\Psi$, then Sym's formula
\BEQ \label{Symformula}
J(\Psi(z,\zquer,\lambda))=-\frac{1}{2H}\left(\pder{F}{\theta}F\inv+\frac{i}{2}F\sigma_3
F\inv\right),\kern1cm\lambda=e^{i\theta},
\EEQ
gives a one parameter family of CMC-immersions in the spinor
representation.
The original CMC-immersion $\Psi$ is reproduced by setting $\lambda=1$.
The whole family $\Psi_\lambda$ of CMC-immersions is usually called
the associated family of $\Psi$.

\newsection{}
The following result was proved in \cite{DoPeWu:1}:

\mytheorem{} {\em 
Let $F:\FD\rightarrow\Lambda_r\LieSU(2)_\sigma$, $\FD\subset\Bcc$, simply
connected, and let $\Psi_\lambda:\FD\rightarrow\threespace$, $\lambda\in
S^1$, be defined by Sym's
formula~\bref{Symformula}. Then the following are equivalent:

1. For each $\lambda\in S^1$ the map 
$\Psi_\lambda:\FD\rightarrow\threespace$ is a CMC-immersion,

2. The $\Lambda_r\Liesu(2)_\sigma$-valued $1$-form 
$\alpha:=F\inv\diff F$ is of the form~\bref{alphaform}.
}

\separate
The construction of such extended frames is the goal of the
DPW method.

To get the so called meromorphic potential from the extended frame we
utilize the Birkhoff splitting:
We define the map
$g_-:\FD\rightarrow\Lambda^-_{r,\ast}\LieSL(2,\Bcc)_\sigma$ by
\BEQ
F(z,\zquer,\lambda)=g_-(z,\zquer,\lambda)g_+(z,\zquer,\lambda)
\EEQ
where $g_+$ takes values in $\Lambda_r^+\LieSL(2,\Bcc)_\sigma$.

The derivative
\BEQ
\xi=g_-\inv\diff g_-
\EEQ
is a one form taking values in
$\Lambda^-_{r,\ast}\Liesl(2,\Bcc)_\sigma$. By a simple calculation 
using the definitions~\bref{UDEF}, \bref{VDEF} of $U(\lambda)$ and
$V(\lambda)$, it is easily shown, that $\xi(z,\lambda)$ is of the form
\BEQ \label{xiform}
\xi(z,\lambda)=\lambda\inv\tmatrix0{f(z)}{\frac{E(z)}{f(z)}}0\diff z,
\EEQ
where $E(z)\diff z^2$ is the Hopf differential and $f(z)$ is a
scalar meromorphic function. The meromorphic potential has poles in
the set $S\subset\FD$ which consists of those points $z\in\FD$, 
for which $F(z,\zquer,\lambda)$ is not in the big cell.
As a function in $\lambda$ it can be extended holomorphically to
$\CPE\setminus\{0\}$.

\newsection{}
Conversely, given a meromorphic potential of the form~\bref{xiform}
one can use the Iwasawa decomposition
\BEQ \label{merotoframe}
g_-(z,\zquer,\lambda)=F(z,\zquer,\lambda)g_+(z,\zquer,\lambda)\inv
\EEQ
at every point in $\FD$ 
to get back a map $F$ from $\FD$ to $\Lambda_r\LieSU(2)_\sigma$. 

A simple calculation then shows, that the $1$-form $\alpha:=F\inv\diff
F$ is of the form \bref{alphaform}.
Therefore, $F(z,\zquer,\lambda)$ defined by Eq.~\bref{merotoframe} is
indeed an extended frame.

It should also be noted, that this construction does not depend on the
chosen radius $r$: Since the meromorphic potential, and therefore also
$g_-(z,\lambda)$ can be continued holomorphically (in $\lambda$) to 
$\CPE\setminus\{0\}$, the extended frame $F(z,\zquer,\lambda)$ can be
continued holomorphically to $\cstar=\Bcc\setminus\{0\}$. Therefore,
the Iwasawa decomposition of $g_-$ at $C_r$ gives always the same
$F$ restricted to the respective circle $C_r$.
This shows that there is only one DPW-method, not a whole family of
``$r$-DPW methods''.

\newsection{} \label{DPWdressing}
Another construction of interest here is the
dressing method. It allows to construct new CMC-immersions from old
ones.

In the loop group formalism the dressing method has a very simple
formulation.

Given an extended frame $F_0$, we take an arbitrary radius $0<r<1$ and
an element $h_+$ of
$\Lambda_r^+\LieSL(2,\Bcc)_\sigma$ and compute the Iwasawa splitting of 
$h_+F_0$ at every point in $\FD$,
\BEQ \label{rdressing}
h_+(\lambda)F_0(z,\zquer,\lambda)=F(z,\zquer,\lambda)p_+(z,\zquer,\lambda),
\EEQ
where $F$ is a new map from $\FD$ to $\Lambda_r\LieSU(2)_\sigma$ and
$p_+$ is a map from $\FD$ to $\Lambda^+_r\LieSL(2,\Bcc)_\sigma$.
Here we choose $p_+$ such that $F$ satisfies $F(0,\lambda)=I$ for all
$\lambda\in S^1$.
Again, an elementary calculation shows, that $F$ is the extended
frame of a (new) CMC-immersion.

If $\xi_0=\lambda\inv\tmatrix0{f_0}{\frac{E_0}{f_0}}0\diff z$ is the
meromorphic potential for $F_0$, then the new meromorphic potential
$\xi$ for $F$ is given by
\BEQ
\xi=\lambda\inv\tmatrix0f{\frac{E}{f}}0\diff z
=p_+\inv\xi_0p_++p_+\inv\diff p_+.
\EEQ
From this it follows, that $E=E_0$, i.e., the Hopf differential
is invariant under dressing.

As in the case of the DPW method it can be easily seen, that for given
$h_+$ this construction does not depend on the chosen radius $r$, {\em
as long as $h_+$ is defined on $C_r$ and can be extended holomorphically
to $I^{(r)}$}. Thus, if we denote the dressing orbit of $F^{0}$ w.r.t.\
$r$ by $\FO_r(F^{0})$, we have a canonical inclusion
\BEQ
\FO_r(F^{0})\subset\FO_{r^\prime}(F^{0}) \kern1cm\mbox{\rm for
$0<r^\prime<r<1$.}
\EEQ

\newsection{} \label{DPWsymmetry} 
Finally, we want to state the transformation properties of an extended
frame under a translation in $\Sym(\Psi)$ as they were derived in
\cite[Section~3.1--3.3]{DoHa:2}. From here on we will set $\FD=\Bcc$.
For notational convenience we will also 
from now on omit the $\zquer$-dependency.

\mytheorem{} {\em
Let $\Psi(z)$ be a CMC-immersion with extended frame 
$F(z,\lambda)$ and associated family $\Psi_\lambda(z)$,
$\Psi(z)=\Psi_1(z)$. 

If $q\in\Sym(\Psi)$ then $F(z,\lambda)$
transforms under the translation by $q$ as
\BEQ \label{Ftrafo}
F(z+q,\lambda)=\chi(q,\lambda)F(z,\lambda)
\EEQ
where $\chi(q,\lambda)$ can be extended holomorphically to $\cstar$. 
For arbitrary $0<r<1$, restriction to $C_r$ gives a group homomorphism
$\chi(\cdot,\lambda):\Sym(\Psi)\rightarrow\Lambda_r\LieSU(2)_\sigma$,
i.e.\
\BEQ \label{chihomo}
\chi(q_1+q_2,\lambda)=\chi(q_1,\lambda)\chi(q_2,\lambda),\kern1cm 
q_1,q_2\in\Sym(\Psi).
\EEQ
}

\Proof By Remark~\ref{intro11}, Eqs.~\bref{Ftrafo} and~\bref{chihomo}
follow for $\lambda\in S^1$ immediately from 
\cite[Lemma~3.1, Theorem~3.3]{DoHa:2}.
Since, by~\bref{Finitial}, $\chi(q,\lambda)=F(q,\lambda)$ we get using
Lemma~\ref{DPWloopgroups}, 
that $\chi(q,\lambda)$ can be extended holomorphically
to $\cstar$. Therefore, restriction to $C_r$, $0<r<1$, gives
$\chi(q,\lambda)\in\Lambda_r\LieSU(2)_\sigma$ and the rest follows by
uniqueness of analytic continuation.
\QED

\separate
As we did already for extended frames, we can, for arbitrary radius
$0<r<1$, view $\chi(q,\lambda)$ as an element of
$\Lambda_r\LieSU(2)_\sigma$ by restricting it to the circle $C_r$.

\newsection{} \label{DPWperiodic}
Using Sym's formula~\bref{Symformula}, we can also easily
compute the proper Euclidean motion $T_q$ in Eq.~\bref{symmetrydef} from the
map $\chi$:
\BEA
J(\Psi(z+q)) & = & -\frac{1}{2H}\bigg(
\pder{\chi(q,\lambda)}{\theta}|_{\theta=0}\cdot\chi(q,1)\inv
+\chi\pder{F(z,\lambda)}{\theta}|_{\theta=0}\cdot
F(z,1)\inv\chi(q,1)\inv\nonumber\\
& & \mbox{}\kern1cm-\frac{i}{2}\chi(q,1)F(z,1)\sigma_3F(z,1)\inv
\chi(q,1)\inv\bigg)\nonumber\\
& = & \chi(q,1)J(\Psi(z))\chi(q,1)\inv
-\frac{1}{2H}\pder{\chi(q,\lambda)}{\theta}|_{\theta=0}\cdot\chi(q,1)\inv,
\kern0.5cm\lambda=e^{i\theta}.\label{Symtrafo}
\EEA
Thus, $T_q$ is given by
\BEQ
J(T_q(z))=J(R_q.z+t_q)
\EEQ
where the rotation $R_q$ is in the spinor representation given by
\BEQ
J(R_q.z)=\chi(q,1)J(z)\chi(q,1)\inv
\EEQ
and the vector $t_q\in\threespace$ is given by
\BEQ
J(t_q)=-\frac{1}{2H}\pder{\chi(q,\lambda)}{\theta}|_{\theta=0}
\cdot\chi(q,1)\inv.
\EEQ
We get the following 

\mytheorem{} {\em
Let $\Psi:\Bcc\rightarrow\threespace$ be a CMC-immersion and let
$q\in\Sym(\Psi)$. Define, for $0<r<1$, 
$\chi(q,\lambda)\in\Lambda_r\LieSU(2)_\sigma$ as above,
then $q\in\Per(\Psi)$ iff 
\BEQ \label{periodicitycond}
\chi(q,1)=\pm I\kern1cm\mbox{\rm and}\kern1cm
\pder{\chi(q,\lambda=e^{i\theta})}{\theta}|_{\theta=0}=0.
\EEQ
}

\Proof Follows from Eq.~\bref{Symtrafo}.
\QED

\separate
From~\bref{Symtrafo} we immediately get the

\mycorollary{} {\em
Let $\Psi:\Bcc\rightarrow\threespace$, $F(z,\lambda)$ and
$\Psi_\lambda$, $\lambda\in S^1$ be as in
Theorem~\ref{DPWsymmetry}. If $F(z,\lambda)$ transforms under the
translation by $q\in\cstar$ as in~\bref{Ftrafo}, then
$q\in\Sym(\Psi_\lambda)$ for all $\lambda\in S^1$.
}

\myremark{}
It actually follows from \cite[Prop.~3.4]{DoHa:2}, that the group
$\Sym(\Psi_\lambda)$ of translational symmetries is the same 
for all elements of the associated family. 
The group $\Per(\Psi_\lambda)$, however, depends crucially on $\lambda$
as we will see in Section~\ref{spaceclosing}.

\newsection{} \label{involutions}
For later use we will introduce the following antiholomorphic involution:
\BEQ \label{taudef}
\tau(\lambda)=\lambdaquer\inv.
\EEQ
Geometrically speaking, $\tau$ is the reflection at the unit
circle. For a map $g(\lambda)$ from a subset of $\CPE$ to 
$\LieSL(2,\Bcc)$ we define
\BEQ
g^\ast(\lambda)=\overline{g(\tau(\lambda))}^\top.
\EEQ
Thus, if $F\in\Lambda_r\LieSL(2,\Bcc)_\sigma$ is defined and holomorphic on
$A^{(r)}$, then $F\in\Lambda_r\LieSU(2)_\sigma$ is equivalent to
\BEQ \label{realitycond}
F^\ast=F\inv.
\EEQ
For a scalar function $f(\lambda)$ we set
\BEQ
f^\ast(\lambda)=\overline{f(\tau(\lambda))}.
\EEQ
If $f$ is defined and holomorphic on a 
$\tau$-invariant neighbourhood of $S^1$, then
$f$ is real on $S^1$ iff $f^\ast=f$.

\section{Translational symmetries} \label{transsym} \message{[transsym]}
In this chapter we will give necessary and
sufficient conditions for a surface to be invariant under a given
translation $z\mapsto z+q$, $q\in\cstar$. The conditions will be
formulated in terms of certain functions that will be introduced in
the next sections.

The calculations here are similar to those in
\cite[Section~3.7]{DoHa:2}, with the exception, that we use here the
$r$-loop formalism recalled in the previous section
from~\cite{DoWu:1}. In \cite{DoHa:2} these
calculations were used to prove, that there are no translationally symmetric 
surfaces in the standard
dressing orbit of the cylinder. However, by the result of
\cite[Corollary~5.3]{DoWu:1} (see also \cite{BuPe:1}), every CMC-torus is
contained in the $r$-dressing orbit of the cylinder.

\newsection{} \label{cylinderorbit}
Let $F(z,\lambda)$ be defined by the $r$-dressing action of
$h_+\in\Lambda^+_r\LieSL(2,\Bcc)_\sigma$
on the extended frame of the (standard) cylinder, i.e.,
\BEQ \label{Fdressingdef}
h_+(\lambda)e^{(\lambda\inv z-\lambda\zquer)A}=F(z,\lambda)\cdot
p_+(z,\lambda),\kern5mm
p_+(z,\lambda)\in\Lambda^+_r\LieSL(2,\Bcc)_\sigma,\kern5mm A=\tmatrix0110.
\EEQ
We note that $e^{(\lambda\inv z-\lambda\zquer)A}$ is the extended
frame associated with the meromorphic potential\
\BEQ
\xi=\lambda\inv\tmatrix0110\diff z
\EEQ
and describes the (standard) cylinder. We also note that all functions
of $z$ are defined on $\Bcc$.

The frame $F(z,\lambda)$ transforms under the translation $z\mapsto
z+q$, $q\in\cstar$, as
\BEQ
F(z+q,\lambda)=Q(q,\lambda)F(z,\lambda)r_+(q,z,\lambda),
\EEQ
where
\BEQ
Q(q,\lambda)=h_+e^{(\lambda\inv q-\lambda\qquer)A}h_+\inv
\EEQ
and
\BEQ
r_+(q,z,\lambda)=p_+(z,\lambda)p_+(z+q,\lambda)\inv.
\EEQ

\newsection{} \label{chiform}
By Theorem~\ref{DPWsymmetry} and Corollary~\ref{DPWperiodic},
we have that $q\in\Sym(\Psi)$, iff
\BEQ \label{Qchirel}
Q(q,\lambda)F(z,\lambda)r_+(q,z,\lambda)=\chi(q,\lambda)F(z,\lambda),
\EEQ
where
\BEQ
\chi(q,\lambda)=F(q,\lambda)
\EEQ
is unitary on $S^1$.
From Section~\ref{involutions} we thus know
\BEQ
\chi^\ast(q,\lambda)=\chi(q,\lambda)\inv \kern1cm\mbox{\rm for
$\lambda\in\cstar$.}
\EEQ
We will derive further conditions on the matrix function $\chi$.

The initial condition~\bref{Finitial} together with Eq.~\bref{Qchirel}
implies 
\BEQ \label{chiQrel}
\chi=Q\cdot R_+,
\EEQ
with
\BEQ
R_+(q,\lambda)=r_+(q,0,\lambda)\in\Lambda_r^+\LieSL(2,\Bcc)_\sigma.
\EEQ
Thus, $F(z,\lambda)$ is invariant under the $r$-dressing transformation
with $R_+$,
\BEQ
R_+(q,\lambda)F(z,\lambda)=F(z,\lambda)r_+(q,z,\lambda).
\EEQ
Substituting~\bref{Fdressingdef} into this equation and rearranging
terms yields
\BEQ \label{dreizweisechs}
e^{-\lambda\inv zA}h_+\inv R_+h_+e^{\lambda\inv zA}=
e^{-\lambda\zquer A}p_+\inv r_+p_+e^{\lambda\zquer A}.
\EEQ
Abbreviating $w_+(\lambda)=h_+\inv R_+h_+$ and
$V_+(z,\lambda)=e^{-\lambda\zquer A}p_+\inv r_+p_+e^{\lambda\zquer
A}$, this is
\BEQ \label{UpVp}
e^{-\lambda\inv zA}w_+(\lambda)e^{\lambda\inv zA}=V_+(z,\lambda).
\EEQ
The l.h.s.\ of~\bref{UpVp} is
\BEQ
Z=(\exp(-\lambda\inv z\ad A))w_+(\lambda).
\EEQ
Thus, \bref{UpVp} says that $Z$ has a Fourier expansion, where only
non-negative powers of $\lambda$ occur. Expanding
$w_+(\lambda)=w_0+w_1\lambda+w_2\lambda^2+\ldots$, altogether we have
\BEQ \label{expansion}
(\exp(-\lambda\inv z\ad A))(w_0+w_1\lambda+w_2\lambda^2+\ldots)
\kern5mm\mbox{\rm does not contain any negative powers of $\lambda$.}
\EEQ
Expanding \bref{expansion} as a function of $z$, a straightforward
induction on the power $N$ of $z$ shows
\BEQ \label{Awpcommute}
[A,w_+(\lambda)]=0.
\EEQ
This implies
\BEQ
V_+(z,\lambda)=V_+(0,\lambda)=w_+(\lambda).
\EEQ
By the definition of $w_+$ we have
\BEQ \label{Rpwp}
R_+=h_+w_+h_+\inv.
\EEQ
Substituting~\bref{Rpwp} into~\bref{chiQrel} finally gives
\BEQ \label{chiwp}
\chi(q,\lambda)=h_+(\lambda)e^{(\lambda\inv q-\lambda\qquer)A}w_+(\lambda)
h_+\inv(\lambda).
\EEQ

\newsection{} \label{Adigression}
This section is a ``digression on $A$''. We consider the matrix
\BEQ \label{Ddef}
D=\frac{1}{\sqrt{2}}\tmatrix11{-1}1
\EEQ
and verify
\BEQ \label{Adiagonal}
DAD\inv=\sigma_3.
\EEQ
Conjugating~\bref{dreizweisechs} by $D$ yields
\BEQ
e^{-\lambda\inv z\sigma_3}\tw_+e^{\lambda\inv
z\sigma_3}\in\Lambda_r^+\LieSL(2,\Bcc)_\sigma,
\EEQ
where $\tw_+=Dw_+D\inv$. Thus, $\tw_+$ and
\BEQ \label{twpinplus}
e^{-\lambda\inv z\sigma_3}\tw_+e^{\lambda\inv z\sigma_3}
\EEQ
are defined at $\lambda=0$.
If we set
\BEQ
\tw_+=\tmatrix{\tw_a}{\tw_b}{\tw_c}{\tw_d},
\EEQ
then the matrix in \bref{twpinplus} has the form
\BEQ
\tmatrix{\tw_a}{e^{-2\lambda\inv z}\tw_b}{e^{2\lambda\inv z}\tw_c}{\tw_d}.
\EEQ
Therefore, if the off-diagonal entries of $\tw_+$ do not vanish
identically, then, for all $z\in\FD$ where $\tw_+$ is defined,
the off-diagonal entries of the matrix in
Eq.~\bref{twpinplus} have an essential singularity at
$\lambda=0$. Therefore, $\tw_b=\tw_c=0$, whence $\tw_+$ 
commutes with $\sigma_3$ and $w_+=D\inv\tw_+D$ commutes with $A$,
reproducing~\bref{Awpcommute}.

The diagonalization~\bref{Adiagonal} is also useful to represent $w_+$
as an exponential. To this end we note, that $\tw_+=Dw_+D\inv$ is a
diagonal matrix $\tw_+=\diag(s_+,s_+\inv)$, with
$s_+(\lambda)\in\cstar$ for $\lambda\in \overline{I^{(r)}}$. 
The map $\overline{I^{(r)}}\rightarrow\cstar$ can be factored through
the universal cover $\Bcc$ of $\cstar$ by standard topological
reasons~\cite[Satz~4.9]{Fo:1}, 
where the covering map $\Bcc\rightarrow\cstar$ is
$\mu\mapsto e^\mu$. Thus, $s_+$ can be written as $s_+=\exp(f_+)$,
where
\BEQ \label{fpholext}
f_+:C_r\rightarrow\Bcc\kern5mm\mbox{\rm is odd and has a holomorphic extension
to $I^{(r)}$.}
\EEQ
Undoing the diagonalization of $w_+$ we obtain
\BEQ \label{wpexp}
w_+=e^{f_+A}.
\EEQ
Finally, in preparation of applications in the rest of this paper, we
discuss the centralizer $C(A)$ of $A$ in
$\Lambda_r\LieSL(2,\Bcc)_\sigma$ in some more detail.
First we note
\BEQ \label{CA}
C(A)=\{\alpha I+\beta A;\alpha,\beta\in\FA_r,\alpha^2-\beta^2=1\}.
\EEQ
Here $\FA_r$ denotes the algebra of functions on $C_r$ that defines
the topology of the various loop groups considered in this paper.

The identity~\bref{CA} is equivalent with the obvious fact, that $A$
is a regular semisimple matrix.

Conjugation with $D$ transforms $A$ into $\sigma_3$ and $C(A)$ into
the centralizer of $\sigma_3$, i.e., into diagonal matrices. Thus,
\BEQ \label{Hdiag}
D(\alpha I+\beta A)D\inv=\diag(\alpha+\beta,\alpha-\beta).
\EEQ
This shows in particular that the eigenvalues of $\alpha I+\beta A\in
C(A)$ are $\alpha+\beta$ and $\alpha-\beta$, and also
$\alpha-\beta=(\alpha+\beta)\inv$, a trivial consequence of
$\alpha^2-\beta^2=1$.

\newsection{} \label{Hform}
Next we consider the matrix
\BEQ \label{Hdef}
H=e^{(\lambda\inv q-\lambda\qquer)A}w_+
\EEQ
occuring in the description~\bref{chiwp} of
$\chi$. From~\bref{Awpcommute} it follows, that $H$ commutes with $A$,
i.e.~\bref{CA} implies
\BEQ \label{alphabetadef}
H=\alpha I+\beta A,
\EEQ
\BEQ \label{a2b2}
\alpha^2-\beta^2=1.
\EEQ
This implies
\BEQ \label{Hinv}
H\inv=\alpha I-\beta A.
\EEQ
Since $q\neq0$, we get from~\bref{Hdef}, that $H$ has an essential
singularity at $\lambda=0$. From this it follows that
\BEQ \label{betanontrivial}
\beta\not\equiv0.
\EEQ
By the twisting condition we also know, that 
\BEQ \label{alphabetaparity}
\mbox{\rm $\alpha$ is even in $\lambda$ and $\beta$ is odd in
$\lambda$.}
\EEQ
Since
\BEQ \label{chiHform}
\chi=h_+Hh_+\inv
\EEQ
is holomorphic for $\lambda\in\cstar$, taking the trace of $\chi$ and
$\chi^2$ yields
\BEQ \label{alphabetaext}
\mbox{\rm $\alpha$ and $\beta^2$ have holomorphic extensions to
$\cstar$.}
\EEQ
We note, that \bref{Hdef} and~\bref{alphabetadef} imply
\BEQ
\mbox{\rm $\beta$ has a holomorphic extension to $0<|\lambda|<r$.}
\EEQ

\myremark{}
In general, $\beta$ will not have a holomorphic extension to
$\cstar$. However, if $\beta$ does have a holomorphic extension to
$\cstar$, then $H$ is also defined on $\cstar$.

Finally, we know from~\bref{wpexp} that $w_+=\exp(f_+A)$. As a
consequence, altogether we have
\BEQ \label{Heq}
H=e^{pA},\kern1cm p=\lambda\inv q-\lambda\qquer+f_+.
\EEQ
In particular, with $f_+$ also $p$ is odd in $\lambda$.

\newsection{} \label{Sdef}
It will be convenient to consider the matrices
\BEQ \label{Sdefeq}
S=h_+Ah_+\inv=\tmatrix abcd,
\EEQ
\BEQ
\hS=\beta S=\tmatrix\ha\hb\hc\hd.
\EEQ
In view of~\bref{alphabetadef} and~\bref{chiHform} we thus have
\BEQ
\chi=\alpha I+\beta S=\alpha I+\hS
\EEQ
and~\bref{Hinv} implies
\BEQ
\chi\inv=\alpha I-\beta S=\alpha I-\hS.
\EEQ
Since, by Theorem~\ref{DPWsymmetry}, 
$\chi$ is holomorphic on $\cstar$, we have 
\BEQ \label{hahbhcext}
\mbox{\rm $\ha$, $\hb$, $\hc$, $\hd$ are holomorphic on $\cstar$.}
\EEQ
Clearly, we have $\tr S=0$ and $\tr\hS=0$, whence
\BEQ \label{deq}
d=-a,\kern3cm \hd=-\ha.
\EEQ
Also, in view of~\bref{alphabetaparity}, the twisting condition for
$\Lambda_r\LieSL(2,\Bcc)_\sigma$ implies
\BEQ \label{abparity}
\mbox{$\ha,b,c$ are even in $\lambda$, $a,\hb,\hc$ are odd in
$\lambda$.}
\EEQ
Two important conditions still need to be evaluated: $S^2=I$ and
$\chi|_{S^1}$ unitary.
The first condition is in view of~\bref{deq} equivalent with
\BEQ \label{abcrel}
a^2+bc=1.
\EEQ
Since $\hS=\beta S$, we also have
\BEQ \label{hahbbeta}
\ha^2+\hb\hc=\beta^2.
\EEQ
The second condition is in view of Section~\ref{involutions}
and~\bref{Hinv} equivalent with
\BEQ \label{alphareality}
\mbox{\rm $\alpha$ is defined on $\cstar$ and $\alpha=\alpha^\ast$,}
\EEQ
\BEQ \label{hSreality}
\mbox{\rm $\hS$ is defined on $\cstar$ and $\hS^\ast=-\hS$.}
\EEQ
In particular,~\bref{hSreality} is equivalent with
\BEQ \label{hahbreal}
\ha^\ast=-\ha,\kern1cm\hb^\ast=-\hc\kern1cm\mbox{\rm for $\lambda\in\cstar$.}
\EEQ
Next we consider the squares of $\alpha,\beta,a,b,c$ and
$\ha,\hb,\hc$. First we note
\BEQ \label{a2b2real}
\mbox{\rm $\alpha^2$ and $\beta^2$ are holomorphic on $\cstar$ and
real on $S^1$.}
\EEQ
This follows for $\alpha^2$ from~\bref{alphareality} and for $\beta^2$
from~\bref{a2b2}. Next,~\bref{hahbreal} implies
\BEQ
\mbox{\rm $\ha^2$ is holomorphic on $\cstar$ and non-positive on
$S^1$.}
\EEQ
Substituting this and~\bref{hahbreal} into~\bref{hahbbeta} gives
\BEQ \label{b2nonpositive}
\mbox{\rm $\beta^2$ is non-positive on $S^1$.}
\EEQ
Since $\ha=\beta a$, we know $\ha^2=\beta^2a^2$. In particular,
$a^2=\frac{\ha^2}{\beta^2}$ is meromorphic on $\cstar$ and real on
$S^1$. Since $a^2$ is by definition also holomorphic at $\lambda=0$, it
follows that $a^2$ is also finite at $\lambda=\infty$. Thus
\BEQ \label{a2rational}
\mbox{\rm $a^2$ is a rational function, real and non-negative on $S^1$, and
finite at $0$ and $\infty$.}
\EEQ
For $b^2$ and $c^2$ one argues
similarly. E.g.~$b^2=\frac{\hb^2}{\beta^2}$ is clearly meromorphic on
$\cstar$ and is also, by the definition of $b$, holomorphic at
$\lambda=0$. From~\bref{hahbreal} we obtain that
$(b^2)^\ast=\frac{(\hb^2)^\ast}{\beta^2}=\frac{\hc^2}{\beta^2}=c^2$ is
also holomorphic at $\lambda=0$. This shows, that $b^2$ is meromorphic
on $\CPE$ and thus rational. Altogether we have shown
\BEQ \label{b2c2rational}
\mbox{\rm $b^2$ and $c^2$ are rational and finite at $0$ and
$\infty$.}
\EEQ
Moreover,
\BEQ \label{b2c2rel}
(b^2)^\ast=c^2.
\EEQ
Next, from~\bref{Sdefeq} we see that $a(\lambda=0)=0$ and
$b(\lambda=0)=c(\lambda=0)\inv$. Since $a$ is odd in $\lambda$, we obtain
\BEQ \label{azero}
\mbox{\rm $a^2$ has a zero of order $2(2n-1)$ for some $n>0$ at
$\lambda=0$,}
\EEQ
\BEQ \label{bc1}
b(\lambda=0)c(\lambda=0)=1.
\EEQ
We also note that the relations 
\BEQ
\ha^2=\beta^2a^2,\kern1cm\hb^2=\beta^2b^2,\kern1cm\hc^2=\beta^2c^2
\EEQ
show that $a^2$, $b^2$, and $c^2$ can have poles only where 
$\beta^2$ has a zero.

Finally, from~\bref{hSreality} we obtain $(\beta b)^\ast=-(\beta
c)$. Hence~\bref{abcrel} implies
\BEQ
\beta^2=\beta^2a^2+\beta b\cdot\beta c=\beta^2a^2-(\beta b)(\beta
b)^\ast.
\EEQ
Therefore, on $S^1$ we obtain $\beta^2(a^2-1)=|\beta b|^2$. Since
$\beta^2$ is non-positive on $S^1$ by~\bref{b2nonpositive}, and
$\beta^2\not\equiv0$ by~\bref{betanontrivial}, we have
$a^2-1\leq0$. Thus,
\BEQ \label{aonS}
0\leq a^2(\lambda)\leq 1 \kern1cm\mbox{\rm for $\lambda\in S^1$}.
\EEQ

\newsection{} \label{lastnecc}
In the last section we considered the matrix $S=h_+Ah_+\inv=\tmatrix
abcd$ and we listed properties of $a,b,c,d$. In the rest of this paper
we will characterize $q\in\Sym(\Psi)$ in terms of $a,b,c,d$. Let us
denote by $\FA_r^+$, $0<r<1$, the subalgebra of those functions in 
$\FA_r$ which can be extended holomorphically to $I^{(r)}$. Then
$a,b,c,d\in\FA_r^+$. To make
sure, that there is also an $h_+$ satisfying~\bref{Sdefeq},
producing $F$ ---and thus $\Psi$--- for which $q\in\Sym(\Psi)$, we
prove

\mytheorem{} {\em Let $a,b,c,d\in\FA_R^+$, $0<R<1$, where $a,d$ are odd and
$b,c$ are even. Then $S=\tmatrix abcd$ is of the form $S=h_+Ah_+\inv$
for some $0<r\leq R$ and $h_+\in\Lambda_r^+\LieSL(2,\Bcc)_\sigma$ iff
\BEQ
d=-a,
\EEQ
\BEQ \label{abcinth}
a^2+bc=1,
\EEQ
\BEQ
b(\lambda=0)\neq0.
\EEQ
}

\Proof ``$\Rightarrow$'' This direction was already proved in the last
section.

``$\Leftarrow$'' Using the matrix $D$ introduced in~\bref{Ddef}, we obtain
\BEQ
S=h_+Ah_+=h_+D\inv\sigma_3Dh_+\inv=Y\sigma_3Y\inv
\EEQ
where
\BEQ
Y=h_+D\inv.
\EEQ
A straightforward formal calculation shows, that 
\BEQ
Y=\frac{i}{\sqrt{2}}\tmatrix{\sqrt{b}}{\sqrt{b}}{\frac{1-a}{\sqrt{b}}}{
-\frac{1+a}{\sqrt{b}}}\cdot\tmatrix{x(\lambda)}00{x(\lambda)\inv}.
\EEQ
Clearly, $Y$ is (formally) a matrix of eigenvectors of $S$. We get for
$h_+=YD$ again formally
\BEQ
h_+=\frac{i}{2\sqrt{b}}
\tmatrix{(x-x\inv)b}{(x+x\inv)b}{(x+x\inv)-(x-x\inv)a}{(x-x\inv)-(x+x\inv)a}.
\EEQ
For this matrix to be defined in our setup we need that the diagonal
entries are even in $\lambda$ and the off-diagonal entries are odd in
$\lambda$ and that everything is defined for all $0\leq|\lambda|\leq
r$, where $0<r\leq R$ is chosen appropriately.

Since $b$ is even and $b(\lambda=0)\neq0$, $\sqrt{b}$ is defined in a 
sufficiently small neighbourhood of $\lambda=0$ and is even there.
So what we need to achieve is that $(x-x\inv)\sqrt{b}$ is even,
$(x+x\inv)\sqrt{b}$ is odd, $\frac{(x+x\inv)-(x-x\inv)a}{\sqrt{b}}$ is
odd, and $\frac{(x-x\inv)-(x+x\inv)a}{\sqrt{b}}$ is even, and all
these functions are defined on a sufficiently small neighbourhood of
$\lambda=0$. This is equivalent with: $x-x\inv$ is even and $x+x\inv$
is odd and $x$ is invertible in a neighbourhood of $\lambda=0$. It is
straightforward to verify, that $x=-i\sqrt{\frac{1+a}{1-a}}$ satisfies these
conditions. In this case
\BEQ
h_+=\frac{1}{\sqrt{c}}\tmatrix1a0c
\EEQ
and $r$ is determined by $\sqrt{c}$.
\QED

\newsection{} \label{nectheorem}
Let us collect the necessary conditions we have derived in
Sections~\ref{chiform}--\ref{Sdef}:

\mytheorem{} {\em
Let $\Psi:\Bcc\rightarrow\threespace$ be a CMC-immersion with extended 
frame $F(z,\lambda)$, s.t.\ $F(z,\lambda)$ is given by dressing the 
cylinder under the $r$-dressing~\bref{rdressing} with some 
$h_+\in\Lambda_r^+\LieSL(2,\Bcc)_\sigma$. Assume also, that for
$q\in\Bcc$, $q\neq0$, $F(z+q,\lambda)=\chi(\lambda)F(z,\lambda)$,
i.e., $q\in\Sym(\Psi)$.
Define $h_+Ah_+\inv=\tmatrix abcd$. Then $d=-a$
and the functions $a(\lambda)$, 
$b(\lambda)$, and $c(\lambda)$ are in $\FA_r^+$ and satisfy the 
following conditions:
\begin{description}
\item[a)] $a^2$, $b^2$, $c^2$ are rational,
\item[b)] $a$ is odd in $\lambda$, $b$ and $c$ are even in $\lambda$,
\item[c)] $a^2+bc=1$.
\item[d)] $a^2$ is real on $S^1$ and $0\leq a^2\leq 1$ on $S^1$,
\item[e)] $c^2=(b^2)^\ast$.
\end{description}
Furthermore, there exists an odd function $f_+$ in $\FA_r^+$, s.t.\
with $\alpha=\cosh(p)$, $\beta=\sinh(p)$, $p=\lambda\inv
q-\lambda\qquer+f_+$, we have
\begin{description}
\item[a')] $\alpha$ and $\beta^2$ are defined on $\cstar$,
\item[b')] $\alpha$ and $\beta^2$ are real on $S^1$,
\item[c')] $\beta^2$ is non-positive on $S^1$,
\item[d')] the functions $\beta a$, $\beta b$, and $\beta c$
extend holomorphically to $\cstar$. 
\end{description}
The matrix function
$\chi(\lambda)$ is given by $\chi=\alpha I+\beta h_+Ah_+\inv$.
}

\Proof
By the results of the last sections we know for the functions 
$a,b,c,d$, defined by $S=h_+Ah_+\inv=\tmatrix abcd$:
\begin{itemize}
\item $d=-a$: \bref{deq},
\item $a^2$, $b^2$, $c^2$ are rational functions:
\bref{a2rational}, \bref{b2c2rational},
\item $a$ is odd, $b$ and $c$ are even in $\lambda$: \bref{abparity},
\item $a^2+bc=1$: \bref{abcrel},
\item $a^2$ is real on $S^1$ and $0\leq a^2\leq 1$ on $S^1$:
\bref{a2rational}, \bref{aonS},
\item $c^2=(b^2)^\ast$: \bref{b2c2rel}.
\end{itemize}
Since $q\in\Sym(\Psi)$ we have
$F(z+q,\lambda)=\chi(\lambda)F(z,\lambda)$ by Theorem~\ref{DPWsymmetry}.
Moreover, $\chi(\lambda)=h_+e^{pA}h_+\inv$, where
$p=\lambda\inv q-\lambda\qquer+f_+$ by~\bref{chiHform} and~\bref{Heq}. Thus
$\chi=\alpha I+\beta S$, where $\alpha=\cosh(p)$, $\beta=\sinh(p)$ and
\begin{itemize}
\item $\alpha$ and $\beta^2$ are defined on $\cstar$: \bref{alphabetaext},
\item $\alpha$ and $\beta^2$ are real on $S^1$: \bref{alphareality},
\bref{a2b2real},
\item $\beta^2$ is non-positive on $S^1$: \bref{b2nonpositive}.
\end{itemize}
Finally, $\ha=\beta a$, $\hb=\beta b$, and $\hc=\beta c$ 
extend holomorphically to $\cstar$ by~\bref{hahbhcext}.
\QED

\newsection{} \label{sufftheorem}
We have seen in Section~\ref{lastnecc} under what conditions on $a,b,c$, and
$d$ we can find an $h_+\in\Lambda_r^+\LieSL(2,\Bcc)_\sigma$, s.t.\
$\tmatrix abcd=h_+Ah_+\inv$. This then defines a CMC-immersion $\Psi$
via dressing of the trivial solution with $h_+$. In this section we
characterize those $a,b,c,d,\alpha,\beta$ such that a given $q\in\cstar$ is
in $\Sym(\Psi)$.

\mytheorem{} {\em
Let there be given three even rational functions $a^2(\lambda)$,
$b^2(\lambda)$, and $c^2(\lambda)$, which satisfy the following conditions
\begin{description}
\item[a)] $a^2$ is real on $S^1$ and $0\leq a^2\leq 1$ on $S^1$,
\item[b)] $c^2=(b^2)^\ast$,
\item[c)] There exists an $0<r<1$, s.t.\ the restrictions of 
$a^2$, $b^2$, and $c^2$ to $C_r$ are the squares of
functions $a$, $b$, $c$ in $\FA_r^+$,
\item[d)] $a$ is odd, $b$ and $c$ are even in $\lambda$,
\item[e)] $a^2+bc=1$.
\end{description}
In addition, with $r$ as in c), 
we assume that there exists an odd function $f_+$ in
$\FA_{r^\prime}^+$, $0<r^\prime\leq r$, 
such that for $p=\lambda\inv q-\lambda\qquer+f_+$,
$\alpha=\cosh(p)$, $\beta=\sinh(p)$, $q\in\cstar$, we have
\begin{description}
\item[a')] $\alpha$ and $\beta^2$ are defined and holomorphic on $\cstar$,
\item[b')] $\alpha$ and $\beta^2$ are real on $S^1$,
\item[c')] $\beta^2$ is non-positive on $S^1$.
\item[d')] The functions $\beta a$, $\beta b$, and $\beta c$ extend
holomorphically to $\cstar$.
\end{description}
Then there exists $0<r^{\dprime}\leq r^\prime$ and
$h_+\in\Lambda_{r^\dprime}^+\LieSL(2,\Bcc)_\sigma$, such that
$h_+Ah_+\inv=\tmatrix abc{-a}$. Moreover, for the extended frame
$F(z,\lambda)$ defined by $h_+e^{(\lambda\inv
z-\lambda\zquer)A}=F(z,\lambda)p_+(z,\lambda)$, $|\lambda|=r^\dprime$,
we have $F(z+q,\lambda)=\chi(\lambda)F(z,\lambda)$, where $\chi=\alpha
I+\beta h_+Ah_+\inv$ is holomorphic on $\cstar$ and takes values in
$\LieSU(2)$ on $S^1$. In particular, $q\in\Sym(\Psi)$ for the
CMC-immersion $\Psi$ associated with $F(z,\lambda)$ via Sym's formula.
}

\Proof
Assume, that we have functions $a^2,b^2,c^2$ and $f_+,p,\alpha,\beta$,
such that a)--e), a')--d') are satisfied. We first want to apply
Theorem~\ref{lastnecc}. We set $d=-a$ and know~\bref{abcinth} by e). 
Since $a,b,c$ are defined at $\lambda=0$ and since $a$ is odd we have
$a(0)=0$, whence $b(0)\neq0$. Thus, by Theorem~\ref{lastnecc}, there
exists some $0<r^\dprime\leq r^\prime<1$ and some
$h_+\in\Lambda_{r^\dprime}^+\LieSL(2,\Bcc)_\sigma$, s.t.\
$S=h_+Ah_+\inv=\tmatrix abc{-a}$. Next we consider the extended
frame defined by the $r$-dressing $h_+e^{(\lambda\inv
z-\lambda\zquer)A}=F(z,\lambda)p_+(z,\lambda)$ of the cylinder. Recall
that we use in this paper the unique Iwasawa splitting discussed in
Section~\ref{DPWloopgroups}. We also set $\chi=h_+e^{pA}h_+\inv=\alpha
I+\beta S$. From a') and d') it follows, that $\chi$ is defined and
holomorphic on $\cstar$. Since with $f_+$ also $p$ is an odd function
in $\lambda$, $pA\in\Lambda_{r^\dprime}\Liesl(2,\Bcc)_\sigma$, whence
$\chi\in\Lambda_{r^\dprime}\LieSL(2,\Bcc)_\sigma$. As outlined in
Section~\ref{Sdef}, $\chi$ is unitary on $S^1$, iff~\bref{alphareality}
and~\bref{hSreality} are satisfied. But~\bref{alphareality} follows
from a'), b'), and the first part of~\bref{hSreality} is just d'). The
second condition is $(\beta a)^\ast=-(\beta a)$ and $(\beta
c)^\ast=-(\beta b)$. To verify this condition we square $\beta a$,
$\beta b$, and $\beta c$ and obtain $((\beta a)^2)^\ast=(\beta a)^2$
and $((\beta c)^2)^\ast=(\beta b)^2$, since $\beta^2$ and $a^2$ are
real by b'), and $(c^2)^\ast=b^2$ by b). Hence $(\beta a)^\ast=\pm\beta
a$ and $(\beta c)^\ast=\pm\beta b$. If $(\beta a)^\ast=\beta a$, then
$\beta a$ is real on $S^1$ and $\beta^2a^2=(\beta a)^2$ is
non-negative on $S^1$. But a) and c') imply that $\beta^2a^2$ is
non-positive on $S^1$. The only possibility for both conditions to hold
is $\beta a=0$ on $S^1$. But in this case, of course, also $(\beta
a)^\ast=-\beta a$ as desired. For the remaining case we consider e)
and obtain $\beta^2=\beta^2a^2+(\beta b)(\beta c)$. If $(\beta
c)^\ast=+\beta b$, then $\beta^2=\beta^2a^2+|\beta c|^2$ on
$S^1$. Hence $|\beta c|^2=\beta^2(1-a^2)$ implies $\beta^2\equiv0$ or
$1-a^2\leq0$. The first case is not possible in view of the form of
$p$. The second case yields in view of a), that $a^2\equiv1$ on
$S^1$. Hence $a=\pm1$ on $\Bcc$, a contradiction, since $a(0)=0$. Thus
$(\beta c)^\ast=-\beta b$ as required.

Finally, we show $q\in\Sym(\Psi)$. To this end we multiply
$\chi=h_+(\alpha I+\beta A)h_+\inv$ from the right with
$h_+e^{(\lambda\inv z-\lambda\zquer)A}$ and obtain
\BEQ
\chi h_+e^{(\lambda\inv
z-\lambda\zquer)A}=h_+e^{(\lambda\inv(z+q)-\lambda(\overline{z+q}))A}e^{f_+}.
\EEQ
Using the definition of $F_+(z,\lambda)$ thus gives
\BEQ
\chi(\lambda)F(z,\lambda)p_+(z,\lambda)=F(z+q,\lambda)p_+(z+q,\lambda)e^{f_+}.
\EEQ
This shows,
\BEQ
F(z+q,\lambda)=\chi(\lambda)F(z,\lambda)
\EEQ
since the Iwasawa splitting chosen is unique and $\chi$ is unitary.
\QED

\newsection \label{spaceclosing}
It remains to state the closing conditions Eq.~\bref{periodicitycond}
in terms of the functions used in Theorem~\ref{sufftheorem}. These
turn out to be of a very simple form.

\mytheorem{} {\em 
Let $\Psi:\Bcc\rightarrow\threespace$ be a CMC-immersion in the
$r$-dressing orbit of the cylinder with associated family
$\{\Psi_\lambda,\lambda\in S^1\}$. Let
$q\in\Sym(\Psi)=\Sym(\Psi_\lambda)$, $q\neq0$, and
$\chi(\lambda)$ and $\beta^2(\lambda)$ as in Theorem~\ref{nectheorem}.
Restrict $\lambda$ to $S^1$ and denote by $(\cdot)^\prime$ 
differentiation w.r.t.\ $\theta$, where $\lambda=e^{i\theta}$. 
Then for $\lambda_0\in S^1$ we have:
\begin{description}
\item[a)] $\chi(\lambda_0)=\pm I$ iff $\beta^2$ vanishes at $\lambda_0$.
\item[b)] $\chi(\lambda_0)=\pm I$ and $\chi^\prime(\lambda_0)=0$ iff 
$\beta^2$ vanishes at least to fourth order at $\lambda_0$.
\end{description}
}

\Proof
Let us first prove a): Let us define $a,b,c$ and $f_+,p,\alpha,\beta$
as in Theorem~\ref{nectheorem}. Then
$\chi=\alpha I+\beta \tmatrix abc{-a}$. By c) in
Theorem~\ref{nectheorem}, we have
$\det\chi=\alpha^2-\beta^2=1$. Therefore, $\chi(\lambda_0)=\pm I$ implies
$\alpha(\lambda_0)=\frac12\tr(\chi(\lambda_0))=\pm1$ and 
$\beta^2(\lambda_0)=0$. 

Conversely, assume $\beta^2(\lambda_0)=0$ for
$\lambda_0\in S^1$. Then $\alpha(\lambda_0)=\pm1$. 
By d) in Theorem~\ref{nectheorem} we have that $a^2$ has no poles on
$S^1$. By c) and e) in Theorem~\ref{nectheorem}, also $b^2$ and $c^2$
are defined everywhere on $S^1$. This shows, that
\BEQ \label{seccond}
\beta^2(\lambda_0)a^2(\lambda_0)=\beta^2(\lambda_0)b^2(\lambda_0)=
\beta^2(\lambda_0)c^2(\lambda_0)=0.
\EEQ
This together with $\alpha(\lambda_0)=\pm1$ implies $\chi(\lambda_0)=\pm I$.

Now we prove b). Let us assume, that $\chi(\lambda_0)=\pm I$ and
$\chi^\prime(\lambda_0)=0$. By Theorem~\ref{nectheorem} we have
\BEQ \label{chiprime}
\chi^\prime=\tmatrix{\alpha^\prime+(\beta a)^\prime}{(\beta b)^\prime}{
(\beta c)^\prime}{\alpha^\prime-(\beta a)^\prime}.
\EEQ
Therefore, and by the proof of part a),
$(\beta a)^2=\beta^2a^2$, $(\beta b)^2=\beta^2b^2$, and
$(\beta c)^2=\beta^2c^2$ vanish at least to fourth order at
$\lambda_0$. By c), d), and e) in Theorem~\ref{nectheorem} we know,
that on $S^1$ $b^2$, $c^2$, and $a^2=1-|b^2|^2$ are defined and
cannot vanish simultaneously.
Therefore, $\beta^2$ vanishes at least to fourth order at $\lambda_0$.

Conversely, assume $\beta^2$ vanishes at least to fourth order at
$\lambda_0$. Then we have already shown, that $\chi(\lambda_0)=\pm I$,
$\alpha(\lambda_0)=\pm1$. By differentiating $\alpha^2-\beta^2=1$, we get
$\alpha^\prime(\lambda_0)=\pm\beta(\lambda_0)\beta^\prime(\lambda_0)=0$.
Also, $a^2$, $b^2$, and $c^2$ are
holomorphic on $S^1$. Therefore, $\beta^2a^2$, $\beta^2b^2$, $\beta^2
c^2$ all vanish at least to fourth order at $\lambda_0$. This shows,
that $\beta a$, $\beta b$, $\beta c$ vanish at least to second order
at $\lambda_0$. From this and~\bref{chiprime}, the claim follows.
\QED

\mycorollary{} {\em
Let $\Psi$ and $\Psi_\lambda$ be defined as in the theorem above and
let $q\in\Sym(\Psi)$, $q\neq0$. Then the set $\Omega_q=\{\lambda\in
S^1|q\in\Per(\Psi_\lambda)\}$ is discrete.
}

\Proof 
We define $p=\lambda\inv q-\lambda\qquer+f_+$, 
$\alpha=\cosh(p)$ and $\beta=\sinh(p)$ as in 
Theorem~\ref{nectheorem}.
By Theorem~\ref{spaceclosing} we know, that $\beta^2$ vanishes on
the set $\Omega_q$. Therefore, if $\Omega_q$ has an accumulation
point, then the holomorphic function $\beta^2$ vanishes everywhere on
$\cstar$, contradicting~\bref{betanontrivial}.
\QED

\myremark{}
The corollary does not exclude the possibility that there
are surfaces in the generalized dressing orbit of the standard cylinder,
whose associated family consists entirely of periodic
surfaces. The members of the 
associated family of the standard cylinder itself, for example,
are cylinders for all $\lambda\in S^1$.
But the corollary implies, that the
translations under which the surfaces are periodic depend on the
spectral parameter. In the case of the cylinder this can easily be
checked by a direct calculation.

\section{Hyperelliptic curves}
\label{periodicsurfaces} \message{[periodicsurfaces]}
Let $\Psi:\Bcc\rightarrow\threespace$ be a CMC-immersion in the
$r$-dressing orbit of the cylinder with associated family
$\{\Psi_\lambda;\lambda\in S^1\}$. I.e., if $F(z,\lambda)$ is the
extended frame of $\Psi$, then there exists $0<r<1$ and
$h_+\in\Lambda_r^+\LieSL(2,\Bcc)_\sigma$, such that $F(z,\lambda)$ is
given by~\bref{Fdressingdef}. In this chapter we also assume, that
$\Psi$ has a periodic metric, i.e., that there exists 
$q\in\Sym(\Psi)=\Sym(\Psi_\lambda)$, $q\neq0$.
As we have seen in Theorem~\ref{nectheorem}
and~Theorem~\ref{sufftheorem}, we can formulate the periodicity
conditions for the metric of a CMC-immersion in terms of scalar 
functions $a,b,c$ given
by $h_+Ah_+\inv=\tmatrix abc{-a}$.  These functions are holomorphic in
a neighbourhood $I^{(r)}$ of $\lambda=0$.  In this section we will
introduce a nonsingular hyperelliptic curve, on which $a$,
$b$ and $c$ can be viewed as meromorphic functions.  From now on, the
symbols $a,b,c,f_+,p,\alpha,\beta$ refer to the functions introduced
in Theorem~\ref{nectheorem}.

\newsection{} \label{FCintro}
First we define the new variable $\nu=\lambda^2$.
The functions $a^2(\lambda)$, $b^2(\lambda)$, and $c^2(\lambda)$ are
by Theorem~\ref{nectheorem} rational and even. 
We will regard $a^2$, $b^2$, and $c^2$ as rational
functions of $\nu$.
Since by b) in Theorem~\ref{nectheorem}, $a$ is an odd function in
$\lambda$, either $a=0$ (which corresponds to the standard cylinder)
or $a^2(\lambda)$ has a zero of order $2(2n-1)$, $n>0$, at
$\lambda=0$. As a function in $\nu$, $a^2$ has therefore a zero of odd
order $2n-1$ at $\nu=0$.
Let $\nu_1,\ldots,\nu_k$ be the points in the $\nu$-plane where $a^2$ has
a pole of odd order, and let $\nu_{k+1},\ldots,\nu_{k+l}$ be the points in the
$\nu$-plane away from $\nu=0$ where $a^2$ has a zero of odd order.

\mylemma{} {\em None of the points $\nu_1,\ldots,\nu_{k+l}\in\cstar$ 
defined above lies on the unit circle.
}

\Proof
We have to show, that $a^2(\nu)$ has neither a pole nor a zero of
odd order on $S^1$. In fact, it doesn't matter here, if we view $a^2$
as a function of $\lambda$ or $\nu$. By d) in
Theorem~\ref{nectheorem}, we know, that $a^2$ has no poles on
$S^1$. By c) and e) in Theorem~\ref{nectheorem}, we have
$(1-a^2)^2=|b^2|^2$ on $S^1$, which shows that $b^2$ and
$c^2=(b^2)^\ast$ are
defined on $S^1$ and also, that $a^2$, $b^2$, and $c^2$ cannot vanish
simultaneously on $S^1$. If $a^2$ has a zero of odd order at
$\lambda_0\in S^1$, then $b^2(\lambda_0)\neq 0$. By d') of
Theorem~\ref{nectheorem}, $(\beta a)^2$
and $(\beta b)^2$ are squares of holomorphic functions on $\cstar$. Thus, the
function $\beta^2=\frac{(\beta a)^2}{a^2}=\frac{(\beta b)^2}{b^2}$ has
both a zero of odd order and of even order at $\lambda_0$. This
implies $\beta\equiv0$, contradicting~\bref{betanontrivial},
\QED

\myprop{} {\em Let $k,l$ and $\nu_1,\ldots,\nu_{k+l}$ be defined as
above. Then $g=\frac12(k+l)$ is an integer and we can order the points
$\nu_1,\ldots,\nu_{2g}$, such that
\BEQ \label{branchpoints}
\nu_{2n}=\tau(\nu_{2n-1}),\kern1cm |\nu_{2n-1}|<1,\kern1cm
n=1,\ldots,g.
\EEQ
}

\Proof
By Lemma~\ref{FCintro}, we have $|\nu_{2n-1}|\neq1$ for
$n=1,\ldots,k+l$. Since $a^2(\nu)$ is real on $S^1$, we have using
Section~\ref{involutions}, that the
set $B=\{\nu_1,\ldots,\nu_{k+l}\}$ is invariant under the antiholomorphic
involution $\tau:\nu\rightarrow\nuquer\inv$. 
Since $\tau$ has no fixed points off the unit
circle, we get that $B$ consists of pairs $(\nu_n,\tau(\nu_n))$.
This shows, that $k+l$ is even, whence $g=\frac12(k+l)$ is an
integer, and that we can order $\{\nu_1,\ldots,\nu_{2g}\}$, such
that~\bref{branchpoints} holds.
\QED

\separate
In the following we will order the points $\nu_1,\ldots,\nu_{2g}$
always such that~\bref{branchpoints} holds.

Consider the algebraic equation
\BEQ \label{FCdef}
\mu^2=\nu\prod_{k=1}^{2g}(\nu-\nu_k).
\EEQ

\mytheorem{} {\em
The plane affine curve $\tilde{\FC}$
defined by~\bref{FCdef} can be uniquely
extended to a compact Riemann surface $\FC$ of genus $g$. The
meromorphic function $\nu:\tilde{\FC}\rightarrow\Bcc$ extends to a
holomorphic map $\pi:\FC\rightarrow\CPE$ of degree $2$. The
branchpoints of $\pi$ are the roots of $\mu^2$ and the point $\infty$.
}

\Proof
The proof follows immediately from~\cite[Lemma~III.1.7]{Miranda:1}
since $\mu^2$ has odd degree.
\QED

\separate
In other words, \bref{FCdef} is a (nonsingular) hyperelliptic curve,
obtained by compactifying the plane affine curve
$\tilde{\FC}=\FC\setminus\{P_\infty\}$, where
$P_\infty=\pi\inv(\infty)\in\FC$ is a single point.

\myremark{}
1. The rational functions $a^2$, $b^2$ and $c^2$ in
Theorem~\ref{nectheorem} are defined in terms of the extended frame
$F(z,\lambda)$ of $\Psi$. I.e., they are defined for the whole
associated family, not only for $\Psi$. Therefore,
also the hyperelliptic curve $\FC$ is associated to the whole family
$\{\Psi_\lambda;\lambda\in S^1\}$. This corresponds to the fact that
$\Sym(\Psi)=\Sym(\Psi_\lambda)$.

2. In Proposition~\ref{aonFC} it will be shown, that
$a(\nu)=\sqrt{a^2(\nu)}$ can be lifted to a nonconstant meromorphic
function on $\FC$.  It determines the complex structure of $\FC$
uniquely (see \cite[I.1.6]{FaKr:1}, \cite[Satz 8.9]{Fo:1}).

\newsection{} \label{meroonFC}
On $\tilde{\FC}=\FC\setminus\{P_\infty\}$, 
$\nu$ and $\mu$ are holomorphic functions.
Every point $P\in\FC\setminus\{P_\infty\}$ is determined uniquely
by the values of $\nu$ and $\mu$ at $P$ and can thus be identified
with the pair $(\nu(P),\mu(P))$. Formally, we will also write
$P_\infty=(\infty,\infty)$.
In this notation we get $\pi(\nu,\mu)=\nu$.
Let us define the hyperelliptic involution $I$ on $\FC$ by
\BEQ \label{Idef}
I(\nu,\mu)=(\nu,-\mu).
\EEQ
A point on $\FC$ is a branchpoint iff it is mapped by
$\pi$ to a branchpoint on $\CPE$.
Clearly, the branchpoints of $\FC$
are precisely the fixed points of $I$, i.e., the points $(\nu_k,0)$,
$k=1,\ldots,2g$, $P_0=(0,0)$ and $P_\infty$.

Using the representation of $\FC$ given in \cite[Chapter
III]{Miranda:1} it is easy to see that the functions $\mu$ and
$\mu\nu^{-(g+1)}$ are local coordinates on $\FC$ at $P_0$ and
$P_\infty$, respectively. It will also be convenient to use the
coordinate $\lambda=\lambda(\nu,\mu)$ on $\FC$ at $P_0$ given by
$\lambda(\nu,\mu)=\mu/\left(\prod_{j=1}^{2g}(\nu-\nu_j)\right)^\frac12$.
In view of the definition of $\FC$ we clearly have the relation
$\lambda^2=\nu$ near $P_0$. In particular, since $\nu\circ I=\nu$ we
obtain $\lambda\circ I=\pm\lambda$. The fact that $\lambda$ is
injective finally implies $\lambda\circ I=-\lambda$.

Using the change of coordinates
$(\nu,\mu)\rightarrow(\nu\inv,\mu\nu^{-(g+1)})$ and the defining
relation for $\FC$ we can define similarly a coordinate
$\lambda\inv$ near $P_\infty$ such that $(\lambda\inv)^2=\nu\inv$
and $\lambda\inv\circ I=-\lambda\inv$.

Let us investigate the set of meromorphic functions on $\FC$.
By~\cite[Proposition~1.10]{Miranda:1}, every meromorphic function on
$\FC$ can be uniquely written as
\BEQ \label{merofuncFC}
f(\nu,\mu)=f_1(\nu)+f_2(\nu)\mu,
\EEQ
with two rational functions $f_1$, $f_2$.

\myremark{} 
It is clear from the representation~\bref{merofuncFC} 
of meromorphic functions on
$\FC$, that each rational function $f_1(\nu)$ can be lifted to a
meromorphic function on $\FC$ by setting
$f(\nu,\mu)=f_1(\nu)$. Clearly, then $f\circ I=f$ for such a function.

Conversely, if $f:\FC\rightarrow\Bcc$ is meromorphic, then it can be
identified with a rational function $f_1(\nu)$ iff $f_2(\nu)\equiv0$
in~\bref{merofuncFC}, i.e., iff it satisfies
$f\circ I=f$. We will frequently 
use this identification of rational functions in
$\nu$ with $I$-invariant meromorphic functions on $\FC$.

\newsection{} \label{FCinvolutions}
Let us define
\BEQ \label{FSdef}
\FS=\pi\inv(S^1)=\{(\nu,\mu)\in\FC;\nu\in S^1\}.
\EEQ
The set $\FS$ is connected if $g$ is even, and has two connected 
components if $g$ is odd.
Since $\FS$ is contained in $\tilde{\FC}$, we can identify it with a
subset of $\Bcc^2$.
Using the antiholomorphic involution 
$\tau:\nu\rightarrow\nuquer\inv$ defined in Section~\ref{involutions},
we define the map
$\tsigma:\tilde{\FC}\setminus\{P_0\}\rightarrow\tilde{\FC}\setminus\{P_0\}$ by
\BEQ \label{tauonFC}
\tsigma:(\nu,\mu)\longmapsto
(\frac1\nuquer,\left(\frac1\nuquer\right)^{(g+1)}
\left(\prod_{j=1}^{2g}\nu_j\right)^{\frac12}\overline{\mu}).
\EEQ
We will choose the sign of the square root such that the points on
$\FS$ are fixed by $\tsigma$.

The following is well known (see e.g.~\cite{Ja:1}):

\mytheorem{} {\em
The map $\tsigma$ defined by~\bref{tauonFC} can
be extended to an antiholomorphic involution $\FCinvolution$ on $\FC$,
which preserves the points of $\FS\subset\FC$.

Furthermore, $\FCinvolution$ commutes with the hyperelliptic
involution and leaves invariant the set of branchpoints of $\FC$.
}

\Proof
Using Theorem~\ref{FCintro} it is easily checked, that $\tsigma$
defines an antiholomorphic involution on $\tilde{\FC}\setminus\{P_0\}$
By using the coordinates $\mu$ and $\mu\nu^{-(g+1)}$ near $P_0$ and
$P_\infty$ respectively, we get $\tsigma(\mu)=c_0\overline{\mu}$,
$c_0=\left(\prod_{j=1}^{2g}\nu_j\right)^\frac12$, in local coordinates,
whence $\tsigma$ extends to an antiholomorphic involution
$\FCinvolution$ on $\FC$, which maps $P_0$ to $P_\infty$. By the
choice of the square root in~\bref{tauonFC}, $\FCinvolution$ fixes the
points on $\FS$.
$\FCinvolution$ clearly commutes with $I$. If $P$ is a branchpoint of
$\FC$, then $I(P)=P$. Therefore,
$I(\FCinvolution(P))=\FCinvolution(I(P))=\FCinvolution(P)$ and
$\FCinvolution(P)$ is also a branchpoint. 
Thus, $\FCinvolution$ leaves invariant the set of branchpoints of $\FC$.
\QED

\separate
For a scalar function on $\FC$ we also define
\BEQ \label{FCstardef}
f^\ast=\overline{f\circ\FCinvolution}.
\EEQ
Since, by Proposition~\ref{FCinvolutions},
$\FCinvolution$ fixes the points of $\FS$, we get

\mylemma{} {\em
Let $f$ be a meromorphic function defined on a
$\FCinvolution$-invariant subset of $\FC$ which contains $\FS$. 
Then $f$ is real on $\FS$ iff $f^\ast=f$ and
$f$ is purely imaginary on $\FS$ iff $f^\ast=-f$.
}

\newsection{} \label{aonFC}
Let us now investigate the properties of $a^2$ w.r.t.\ $\FC$.

\myprop{} {\em
The rational function $a^2$ defined in Theorem~\ref{nectheorem} is of
the form
\BEQ
a^2(\nu)=f(\nu)^2\mu^2(\nu),
\EEQ
where $f$ is rational and defined at $\nu=0$.
The function $a=f\mu$ is a meromorphic function on 
$\FC$, which satisfies
\BEQ \label{FCaodd}
a\circ I=-a
\EEQ 
and 
\BEQ \label{FCareal}
a^\ast=a.
\EEQ
}

\Proof
By the definition of $\mu^2$, the quotient $\frac{a^2}{\mu^2}$ is
rational and has only poles and zeroes of even order. Therefore, it is
the square of a rational function $f(\nu)$. Since $a^2$ has a zero of
odd order at $\nu=0$, $f(\nu)$ is defined at $\nu=0$.
By~\bref{merofuncFC}, $a=f\mu$ is a meromorphic function on $\FC$.
Since $a^2$ is real and non-negative on $S^1$, the function
$a=\sqrt{a^2}$ takes
real values over $S^1$ on $\FC$, whence, by Lemma~\ref{FCinvolutions},
\bref{FCareal} holds.
Furthermore, $a\circ I=-f\mu=-a$ and~\bref{FCaodd} holds.
\QED

\mylemma{} {\em
Let $\phi$ be a meromorphic function on $\FC$. Then $\phi$ is
antisymmetric w.r.t.\ the hyperelliptic involution 
$I$, i.e.\ $\phi\circ I=-\phi$, iff $\phi(\nu,\mu)=f(\nu)\mu$, where
$f$ is a rational function. Furthermore, in this case the following holds:
\begin{description}
\item[1.] Locally around $P_0$, $\phi$ is an odd
meromorphic function of the local coordinate $\lambda$.
\item[2.] The product $\phi a$ can be identified with a rational function of
$\nu$ on $\CPE$.
\end{description}
}

\Proof
The equivalence statement follows immediately from the
representation~\bref{merofuncFC} of
meromorphic functions on $\FC$.
Now consider $\phi_0=\phi\circ\lambda^{(-1)}$, 
where $\lambda^{(-1)}$ denotes the
inverse of the local coordinate map $\lambda$ around $P_0$.
Noting that~\ref{meroonFC} implies
$\lambda\circ I\circ\lambda^{(-1)}(z)=-z$, we obtain
$\phi_0(-z)=\phi\circ I\circ \lambda^{(-1)}(z)=-\phi_0(z)$, and therefore
1.\ holds.

Finally, by~\bref{FCaodd},
we have that $\phi a$ is invariant under $I$, whence, by
Remark~\ref{meroonFC}, can be identified
with a meromorphic function $(\phi a)(\nu)$ on $\CPE$.
\QED

\newsection{} \label{FCprime}
Let us also define the Riemann surface $\FC^\prime$ on which
$\sqrt{a^2(\lambda)}$ is meromorphic. To be precise (see~\cite[Lemma
III.1.7]{Miranda:1}), let $\FC^\prime$
be the hyperelliptic curve associated with the plane affine curve
$\tilde{\FC}^\prime$ defined by the algebraic equation
\BEQ \label{FCprimedef}
\tmu^2=\prod_{i=1}^{2g}(\lambda-\sqrt{\nu_i})(\lambda+\sqrt{\nu_i}).
\EEQ
The holomorphic function $\lambda$ extends to a holomorphic map
$\pi^\prime:\FC^\prime\rightarrow\CPE$ of degree $2$.

As for $\FC$, we will identify the points of $\FC^\prime$ which don't
lie over $\lambda=\infty$ with pairs
$(\lambda,\tmu)$.  It should be noted, that $\FC^\prime$ has no
branchpoints over $\lambda=0$ and $\lambda=\infty$ since $a^2$ has, as
an even function of $\lambda$, a zero of even order at $0$.
In particular, $(\pi^\prime)\inv(\infty)$ consists of
two different points $P_\infty^{(1)}$ and $P_\infty^{(2)}$.
By $P_0^{(1)}$ and $P_0^{(2)}$ we will denote the two covering points
of $\lambda=0$.
Clearly, $\lambda$ is a local coordinate around $P_0^{(1)}$ and
$P_0^{(2)}$ and $\lambda\inv$ is a local coordinate around
$P_\infty^{(1)}$ and $P_\infty^{(2)}$.

Every meromorphic function $\tf$ on $\FC^\prime$ is of the
form~\cite[Proposition~1.10]{Miranda:1}
\BEQ \label{merofuncFCprime}
\tf(\lambda,\tmu)=\tf_1(\lambda)+\tf_2(\lambda)\tmu
\EEQ
with two rational functions $\tf_1$ and $\tf_2$.

\myremark{ 1} Let $I^\prime(\lambda,\tmu)=(\lambda,-\tmu)$ 
be the hyperelliptic involution on
$\FC^\prime$. As in Remark~\ref{meroonFC}, we will use the
representation~\bref{merofuncFCprime} to identify rational functions
of $\lambda$ with $I^\prime$-invariant meromorphic functions on 
$\FC^\prime$.
For the understanding of the rest of the paper it will be helpful to
investigate the relations between $\FC$ and $\FC^\prime$. In
preparation of this we introduce the notation $\Bcc_\nu$,
$\CnPE$, $\cnstar$ for the $\nu$-plane, the projective
$\nu$-plane and $\cnstar=\Bcc_\nu\setminus\{0\}$, respectively.
Similar notation will be used for the $\lambda$-plane. We define
\BEQ
\rho_0:\ClPE\longrightarrow\CnPE,\kern1cm\rho_0(\lambda)=\lambda^2.
\EEQ

\myremark{ 2} 1. For similarly defined functions on $\Bcc_\lambda$ and
$\clstar$ we will use the same notation $\rho_0$.

2. In the rest of the paper we will carefully state where functions,
$1$-forms etc.\ are defined. Notation, like $\rho_0$, will be used for
several closely related maps, $1$-forms etc.

3. If $g$ is a map defined on $\Bcc_\nu$, then ``$g$ on
$\Bcc_\lambda$'' means $g\circ\rho_0$. Similarly, if $\delta$ is a
$1$-form on $\Bcc_\nu$, ``$\delta$ on $\Bcc_\lambda$'' means
$\rho_0^\ast\delta$, the ``pullback of $\delta$ relative to
$\rho_0$''. Similar conventions will be applied to various other
maps.

\separate
Before investigating the relation between $\FC$ and $\FC^\prime$ we
introduce some more notation
\BEQ
\pi:\FC\rightarrow\CnPE
\EEQ
where $\pi$ is the extension of $(\nu,\mu)\rightarrow\nu$ as stated in
Theorem~\ref{FCintro}. Similarly we will use
\BEQ
\pi^\prime:\FC^\prime\rightarrow\ClPE.
\EEQ
Note that $0$ and $\infty$ are branchpoints for $\pi$ but not for
$\pi^\prime$.

Next we consider the map
\BEQ
\rho:\Bcc_\lambda\times\Bcc_{\tmu}\rightarrow\Bcc_\nu\times\Bcc_\mu,\kern1cm
\rho(\lambda,\tmu)=(\lambda^2,\lambda\tmu).
\EEQ
Also recall the definition of the antiholomorphic involution
$\FCinvolution$ on $\FC$ given by~\bref{tauonFC}. Similarly we have an
antiholomorphic involution $\FCinvolution^\prime$ on $\FC^\prime$
defined on $\tilde{\FC}^\prime\setminus\{P_0^{(1)},P_0^{(2)}\}$ by
\BEQ
\FCinvolution^\prime(\lambda,\tmu)
=(\lambdaquer\inv,\lambdaquer^{-2g}\left(\prod_{j=1}^{2g}\nu_j\right)^\frac12
\tmu).
\EEQ

\mytheorem{} {\em
The map $\rho$, restricted to $\tilde{\FC}^\prime$, extends to a
surjective and holomorphic map $\rho:\FC^\prime\rightarrow\FC$. Moreover,
$\rho\circ\FCinvolution^\prime=\FCinvolution\circ\rho$ and $\rho\circ
I^\prime=I\circ\rho$ holds.
}

\Proof
Let $(\lambda,\tmu)\in\tilde{\FC}^\prime$. Then
$\rho(\lambda,\tmu)=(\lambda^2,\lambda\tmu)$ satisfies
\BEQ
(\lambda\tmu)^2-\lambda^2\prod_{j=1}^{2g}(\lambda^2-\nu_j)
=\lambda^2\left(\tmu^2-\prod_{j=1}^{2g}(\lambda^2-\nu_k)\right)=0.
\EEQ
Conversely, if $(\nu,\mu)\in\tilde{\FC}^\ast$, then
$(\sqrt{\nu},\frac{1}{\sqrt{\nu}\mu})\in\Bcc_\lambda\times\Bcc_{\tmu}$
satisfies 
\BEQ
\left(\frac{1}{\sqrt{\nu}}\right)^2
-\prod_{j=1}^{2g}(\sqrt{\nu}-\sqrt{nu_j})(\sqrt{\nu}+\sqrt{\nu_j})
=\frac{1}{\nu}\mu^2-\prod_{j=1}^{2g}(\nu-\nu_j)=0.
\EEQ
If $(\nu,\mu)\in\tilde{\FC}$, $\nu=0$, then $\mu=0$ and
$\rho(0,\tmu_0)=(0,0)$, where $\tmu_0$ is chosen such that
$(0,\tmu_0)\in\tilde{\FC}^\prime$. Since for every point of
$\tilde{\FC}$ either $\nu$ or $\mu$ is a local coordinate, it is easy
to verify that $\rho:\tilde{\FC}^\prime\rightarrow\tilde{\FC}$ is
holomorphic. Finally, let $P=(\lambda(z),\tmu(z))$ denote the points
in a neighbourhood of one of the points above $\infty$ on
$\FC^\prime$. Then, by \cite[Chapter~III]{Miranda:1}, this point is
described in $\tilde{\FC}^\prime$ by
$P=(\lambda(z)\inv\tmu(z)\lambda(z)^{-g^\prime-1})$ where we have
$4g=2g^\prime+2$, thus $g^\prime=2g-1$. Then
$\rho(P)=(\lambda(z)^{-2},\tmu(z)\lambda(z)^{-g^\prime-2})$. In the
chart around $P_\infty\in\FC$ this is
\BEQ
(\lambda(z)^2,\tmu(z)\lambda(z)^{-g^\prime-2}(\lambda(z)^{-2})^{-(g+1)}
=(\lambda(z)^2,\tmu(z)\lambda(z)).
\EEQ
This shows, that $\rho$ extends to a surjective holomorphic map from
$\FC^\prime$ to $\FC$. The last statement is a straightforward
computation.
\QED

\mycorollary{ 1} {\em The following diagram of surjective holomorphic maps
is commutative.
\BEQ \label{CD}
\begin{CD}
\FC^\prime @>\rho>> \FC\\
@V{\pi^\prime}VV @VV{\pi}V\\
\ClPE @>{\rho_0}>> \CnPE
\end{CD}
\EEQ
}

\separate
As a consequene of the Theorem above we can lift every function on
$\FC$ to a function on $\FC^\prime$ via composition with $\rho$. This
applies in particular to meromorphic functions. Denoting by $\FM(\FC)$
and $\FM(\FC^\prime)$ the fields of meromorphic functions on $\FC$ and
$\FC^\prime$ respectively, we obtain the

\myprop{} {\em a) The map $\hrho:\FM(\FC)\rightarrow\FM(\FC^\prime)$,
$\hrho(f)=f\circ\rho$ is an injective homomorphism of fields.

b) A function $f^\prime\in\FM(\FC^\prime)$ is in the image of $\hrho$
iff
\BEQ \label{fprimecond}
f^\prime(-\lambda,\tmu)=f^\prime(\lambda,-\mu).
\EEQ
}

\Proof
a) straightforward.

b) If $f^\prime$ is in the image of $\hrho$, then
\BEQ
f^\prime(\lambda,\tmu)=(f\circ\rho)(\lambda,\tmu)=f(\lambda^2,\lambda\tmu)
=f_1(\lambda^2)+f_2(\lambda^2)\lambda\tmu,
\EEQ
where we have used~\bref{merofuncFC}. A comparison
with~\bref{merofuncFCprime} now shows
\BEA
f_1^\prime(\lambda) & = & f_1(\lambda^2) \label{f1cond}\\
f_2^\prime(\lambda) & = & f_2(\lambda^2)\lambda, \label{f2cond}
\EEA
from which~\bref{fprimecond} follows.
Conversely, assume that
$f^\prime(\lambda,\mu)=f_1^\prime(\lambda)+f_2^\prime(\lambda)\tmu$
satisfies~\bref{fprimecond}. Then $f_1^\prime$ is even and
$f_2^\prime$ is odd in $\lambda$ and we can find functions $f_1$ and
$f_2$ of $\nu=\lambda^2$ such that~\bref{f1cond} and~\bref{f2cond}
hold. Obviously, $f^\prime=f\circ\rho$ for $f=f_1+f_2\mu\in\FM(\FC)$.
\QED

\mycorollary{ 2} {\em
The function $a^2(\lambda)$ on $\Bcc_\lambda$ is the square of a
meromorphic function on $\FC^\prime$.
}

\Proof
By Proposition~\ref{aonFC}, $a^2(\nu)$ is the square of a meromorphic
function on $\FC$. By the proposition above, we know that
$\ca=a\circ\rho$ is meromorphic on $\FC^\prime$. Moreover, $\ca\circ
I^\prime=a\circ\rho\circ I^\prime=a\circ
I\circ\rho=-a\circ\rho=-\ca$, where we have used~\bref{FCaodd}. This
implies that $(\ca)^2$ descends to the meromorphic function
$(\ca)^2(\lambda)=a^2(\lambda^2)$ on $\Bcc_\lambda$ by
Remark~\ref{meroonFC}.
\QED

\myremark{ 3} 1. The diagram~\ref{CD} naturally induces also diagrams,
where the projective spaces are replaced by $\Bcc_\lambda$ and
$\Bcc_\nu$ or by $\clstar$ and $\cnstar$.

2. Each of the surjective maps of any of the above diagrams induces an
injective map on the level of meromorphic functions via composition.

\newsection{} \label{betaonFC}
Let us define the non-compact Riemann surfaces
$$
\mbox{\rm $\FCstar=\FC\setminus\{P_0,P_\infty\}$ and
${\FC^\prime}^\ast
=\FC^\prime\setminus\{P_0^{(1)},P_0^{(2)},P_\infty^{(1)},P_\infty^{(2)}\}$.}
$$
We already know, that $a$ is meromorphic on $\FC^\prime$ and $\FC$,
and therefore also on $\FCstar$. Now we prove the following important result:

\mytheorem{} {\em 
The functions $\alpha$ and $\beta$ are holomorphic on $\FCstar$.
The functions $b=\sqrt{b^2}$ and $c=\sqrt{c^2}$ 
are meromorphic on $\FC^\prime$ without poles over $0$ and $\infty$.
}

\Proof
We know by a') and d') in Theorem~\ref{nectheorem}, 
that $\beta a$, $\alpha$, and $\beta^2$ are even functions in
$\lambda$, which are defined and holomorphic on
$\cstar$. Therefore, after replacing $\lambda^2$ by $\nu$
Remark~\ref{meroonFC} applies and the functions above can be identified
with $I$-invariant holomorphic functions on $\FCstar$. Since also
$a$ is meromorphic on $\FC$, we have that $\beta=\frac{\beta a}{a}$ is
a meromorphic function on $\FCstar$. Since the square $\beta^2(\nu)$ is
holomorphic on $\cstar$, $\beta$ has no poles on $\FCstar$, which
shows that $\beta$ is holomorphic on $\FCstar$.

Applying~\bref{CD} we pull back $\beta$ from $\FCstar$ to a holomorphic
function (also denoted by $\beta$) on ${\FC^\prime}^\ast$.
Furthermore, $\beta b$ and $\beta c$ are odd functions of $\lambda$,
which are defined and holomorphic on $\cstar$. Thus, $\beta b$, $\beta
c$, $b=\frac{\beta b}{\beta}$ and $c=\frac{\beta c}{\beta}$ are
meromorphic on ${\FC^\prime}^\ast$.
Since, by Theorem~\ref{nectheorem}, $b(\lambda)$ and $c(\lambda)$ 
are in $\FA_r^+$, they
can be continued holomorphically to $\lambda=0$ on $\CPE$.  By e) in
Theorem~\ref{nectheorem}, the same holds for $b$ and $c$ around
$\lambda=\infty$.  Since $\lambda$ is a local coordinate around
$P_0^{(1)}$ and $P_0^{(2)}$ and $\lambda\inv$ is a local coordinate
around $P_\infty^{(1)}$ and $P_\infty^{(2)}$, the functions $b$ and
$c$ can be extended holomorphically to these points, which finishes
the proof.
\QED

\myprop{} {\em
With $(\cdot)^\ast$ defined in Section~\ref{FCinvolutions}, we have
\BEQ \label{461}
\alpha^\ast=\alpha,\kern1cm\beta^\ast=-\beta\kern1cm\mbox{\rm on
$\FC^\ast$ and on ${\FC^\prime}^\ast$},
\EEQ
\BEQ \label{462}
\mbox{\rm $b^\ast=c$ on $\FC^\prime$},
\EEQ
\BEQ \label{463}
\mbox{\rm $\alpha\circ I=\alpha$ and $\beta\circ I=-\beta$ on
$\FCstar$},
\EEQ
\BEQ \label{464}
\mbox{\rm $\alpha\circ I^\prime=\alpha$ and $\beta\circ I^\prime=
-\beta$ on ${\FC^\prime}^\ast$},
\EEQ
\BEQ \label{465}
\mbox{\rm $b\circ I=-b$ and $c\circ I=-c$ on $\FC^\prime$}.
\EEQ
}

\Proof
By b') in Theorem~\ref{nectheorem}, as functions on $\cstar$, $\alpha$
and $\beta^2$ are real on $S^1$, i.e., if we identify $\alpha$ and
$\beta^2$ with holomorphic functions on $\FCstar$, then they are real
on $\FS$ defined in~\bref{FSdef}. Thus, by Lemma~\ref{FCinvolutions}, 
$\alpha^\ast=\alpha$ holds. By c') in
Theorem~\ref{nectheorem}, $\beta^2$ is non-positive on
$S^1$. Therefore, the holomorphic function $\beta$ on $\FCstar$ is on $\FS$
the square root of a non-positive real function. Hence, $\beta$ is
purely imaginary on $\FS$. This implies $\beta^\ast=-\beta$, by
Lemma~\ref{FCinvolutions}.
The statement on ${\FC^\prime}^\ast$ follows from this since $\rho$
intertwines $\FCinvolution$ with $\FCinvolution^\prime$ by
Theorem~\ref{FCprime}.
To verify \bref{462} we use again the representation $b=\frac{\beta
b}{\beta}$ on ${\FC^\prime}^\ast$ and~\bref{hahbreal} and obtain
$b^\ast=\frac{(\beta b)^\ast}{\beta^\ast}=\frac{\beta
c}{\beta}=c$. The first relation in~\bref{463} is clear, since $\alpha$
is defined on $\cnstar$. For the second we use
$\beta=\frac{\beta a}{a}$ and obtain $\beta\circ I=\frac{(\beta
a)\circ I}{a\circ I}=-\frac{\beta a}{a}=-\beta$, since $\beta a$ is
defined on $\cnstar$ and~\bref{FCaodd} holds. \bref{464} follows
from~\bref{463} since $\rho$ intertwines $I$ with $I^\prime$. Finally,
$b\circ I=\frac{(\beta b)\circ I}{\beta\circ I}=-\frac{\beta
b}{\beta}=-b$. and similarly $c\circ I=-c$.
\QED

\newsection{} \label{eigenbundle}
Let us pause to collect some of the implications of Theorem~\ref{betaonFC}.
By~\bref{Adiagonal}, the matrix 
$S=h_+Ah_+\inv=\tmatrix abc{-a}$ can be written as
\BEQ
S=h_+D\inv\tmatrix100{-1}Dh_+\inv.
\EEQ
The entries of $S$ are all meromorphic functions on $\FC^\prime$.
By linear algebra over the field $\FM(\FC^\prime)$ of meromorphic
functions on $\FC^\prime$, we get that there exist nonzero elements $x_\pm$
of $\FM(\FC^\prime)^2$, written as column vectors,
such that $Sx_\pm=\pm x_\pm$. These vectors are unique up to
multiplication by elements of $\FM(\FC^\prime)$. They define a
two-dimensional vector-bundle over $\FC^\prime$, the ``eigenbundle'' of
$S$ over $\FC^\prime$. If we fix them by requiring that
\BEQ
(x_+(\lambda),x_-(\lambda))=h_+(\lambda)D\inv \kern5mm\mbox{\rm for
$|\lambda|\leq r$},
\EEQ
then they define a lift of $h_+$ to a matrix function with
entries meromorphic on $\FC^\prime$.

Since, by~\bref{chiHform} and~\bref{Hdiag},
the matrix $\chi=\alpha I+\beta h_+Ah_+\inv$ can be written in the form 
\BEQ
\chi=h_+D\inv\tmatrix{\alpha+\beta}00{\alpha-\beta}Dh_+\inv,
\EEQ
we see, that the branchpoints of $\FC^\prime$ are points, where
the eigenvalues of $\chi$ coincide, i.e., where $\chi$ fails to be
regular semisimple.

\newsection{} \label{higherflows}
The next two sections are a digression on finite type surfaces.

Let $\Psi:\Bcc\rightarrow\threespace$ be a CMC-immersion with extended
frame $F(z,\lambda)$ which is generated by $r$-dressing the cylinder with
$h_+\in\Lambda_r^+\LieSL(2,\Bcc)_\sigma$.
Let $\FZ$ be the abelian Lie-subalgebra
of all $\zeta\in\Lambda_r\Liesl(2,\Bcc)_\sigma$,
such that
\begin{itemize}
\item $[\zeta(\lambda),A]=0$ for all $\lambda\in C_r$ and
\item $\zeta$ can be extended to a meromorphic function on $I^{(r)}$ 
which has precisely one pole at $\lambda=0$.
\end{itemize}

\myremark{ 1}
By the discussion in Section~\ref{Adigression} and the
tracelessness of the elements of $\Lambda_r\Liesl(2,\Bcc)_\sigma$,
we can write each $\zeta\in\FZ$ uniquely as 
\BEQ \label{FZrep}
\zeta(\lambda)=\phi(\lambda)A,
\EEQ
where $\phi\in\FA_r$ is odd in $\lambda$, and has precisely
one pole at $\lambda=0$.

\separate
Let $F:\Bcc\rightarrow\Lambda_r\LieSU(2)_\sigma$ be an extended frame in
the $r$-dressing orbit of the cylinder. Thus, there is
$h_+\in\Lambda_r^+\LieSL(2,\Bcc)_\sigma$ and
$p_+:\Bcc\rightarrow\Lambda_r^+\LieSL(2,\Bcc)_\sigma$, such that
\BEQ
F=h_+e^{(\lambda\inv z-\lambda\zquer)A}p_+(z,\lambda)\inv.
\EEQ
For each $\zeta\in\FZ$ we define
\BEQ \label{hpaction}
h_+\#\zeta=U\inv h_+e^{\zeta(\lambda)},
\EEQ
where $U\in\Lambda_r\LieSU(2)_\sigma$ 
is determined by the Iwasawa decomposition~\bref{Iwasawa}, such that
$h_+\#\zeta\in\Lambda_r^+\LieSL(2,\Bcc)_\sigma$.
The action of $\zeta$ on $h_+$ descends
via~\bref{Fdressingdef} to an action on the frame $F$:
\BEQ \label{Zaction}
F\#\zeta
=(h_+\#\zeta)e^{(\lambda\inv z-\lambda\zquer)A}\tp_+(z,\lambda)\inv
=U\inv h_+e^{(\lambda\inv z-\lambda\zquer)A+\zeta(\lambda)}
\tp_+(z,\lambda)\inv,
\EEQ
where $\tp_+:\Bcc\rightarrow\Lambda_r^+\LieSL(2,\Bcc)_\sigma$
is chosen such that $F\#\zeta\in\Lambda_r\LieSU(2)_\sigma$ and
$(F\#\zeta)(0,\lambda)=I$.

Our definition of $\FZ$ and the associated action on the $r$-dressing
orbit of the standard cylinder are obtained from the definitions in
\cite[Section~4]{BuPe:1} for the special case of harmonic maps into
$\LieSU(2)/\LieU(1)$.
By~\cite[Prop.~4.1]{BuPe:1}, we have

\mylemma{} {\em
1. Eqs.~\bref{hpaction} and \bref{Zaction} define an action of $\FZ$ on the
$r$-dressing orbit of the cylinder which is compatible with the
linear structure of $\FZ$. In particular for
$\zeta,\zeta^\prime\in\FZ$, 
\BEQ \label{commute}
F\#(\zeta+\zeta^\prime)=(F\#\zeta)\#\zeta^\prime=(F\#\zeta^\prime)\#\zeta.
\EEQ

2. If $\zeta,\zeta^\prime\in\FZ$, such that $\zeta-\zeta^\prime$ can be
extended holomorphically to $I^{(r)}$, then $F\#\zeta=F\#\zeta^\prime$.
}

\separate
In other words, each $\zeta\in\FZ$ generates a ``higher flow'' on the set of 
CMC-immersions in the $r$-dressing orbit of the cylinder by
\BEQ \label{flowfromaction}
F(z,t,\lambda)=F(z,\lambda)\#(t\zeta),\kern1cm t\in\Brr.
\EEQ
For each $t\in\Brr$, Sym's formula gives a new
associated family of CMC-immersions with extended frame $F(z,t,\lambda)$.
These higher flows commute by~\bref{commute}.

\myremark{ 2}
1. The extended frames $F\#\zeta$ and $F(z,t,\lambda)$ all
satisfy the normalization condition~\bref{Finitial}.

2. Using~\bref{UDEF}--\bref{VF} and the equivalence of~\bref{ZCC}
and~\bref{GaussCodazzi1},\bref{GaussCodazzi2}, we see, 
that~\bref{hpaction}, \bref{Zaction} and~\bref{flowfromaction} also
defines an action on the set of real valued solutions $u(z,\zquer)$
of~\bref{GaussCodazzi1}. Since $E\equiv1$ for all surfaces in the
$r$-dressing orbit of the cylinder, we get by setting $H=-2$,
that this defines a hierarchy of commuting flows of the
integrable $\sinh$-Gordon equation.

\newsection{} \label{finitetype}
We say, that the flow defined by $\zeta\in\FZ$ acts trivial on an
extended frame $F$ in the dressing orbit of the cylinder, iff
\BEQ \label{trivialflow}
F(z,t,\lambda)=U_0(t)F(z,\lambda)U_0(t)\inv\kern1cm\mbox{\rm 
for all $t\in\Brr$,}
\EEQ
where $U_0(t)\in\LieU(1)$ is a $\lambda$-independent unitary matrix.
By~\cite[Corollary~4.1]{DoHa:2}, this is equivalent to the fact that
the associated families $\Psi(z,t,\lambda)$ and $\Psi(z,\lambda)$
defined by $F(z,t,\lambda)$ and $F(z,\lambda)=F(z,0,\lambda)$, differ
only by a proper Euclidean motion in $\threespace$.
Clearly, the set $\FZ^\prime\subset\FZ$ of elements which
generate trivial flows, is a linear subspace of $\FZ$.

\mydefinition{}
A CMC-surface in the $r$-dressing orbit of the cylinder is of
{\em finite type} iff the subspace $\FZ^\prime$ of trivial flows has
finite codimension in $\FZ$. 

\separate
Let $\zeta\in\FZ^\prime$. Then, by~\cite[Corollary~2.6]{DoHa:2},
all surfaces generated by the $\zeta$-flow have the same metric. Therefore,
using Remark~2 in Section~\ref{higherflows}, we see, that $\zeta$
generates a trivial flow of the $\sinh$-Gordon
equation, i.e., $u(z,\zquer,t)=u(z,\zquer)$ for all $t\in\Brr$. 
In~\cite[Theorem~4.2]{BuPe:1} it was shown, that the conformal
factor $u(z,\zquer)$ associated to a CMC-immersion $\Psi$
is a finite type solution of the $\sinh$-Gordon-equation (for a
definition see~\cite{PiSt:1} or~\cite{DoWu:1}) iff $\Psi$ is of
finite type in the sense of Definition~\ref{finitetype}.

In~\cite{PiSt:1} it was shown, that for a surface
with doubly periodic metric, in particular for CMC-tori,
the conformal factor $u$ is a finite type solution of the
$\sinh$-Gordon equation. For more general CMC-surfaces with periodic
metric, e.g., for CMC-surfaces which are topological
cylinders, this is in general not true.
However, for surfaces in the $r$-dressing orbit of the
cylinder, we have the remarkable

\mytheorem{} {\em
Every CMC-surface with periodic metric in the $r$-dressing orbit of
the cylinder is of finite type.
}

\Proof
Consider a CMC-surface $\Psi:\Bcc\rightarrow\threespace$ with periodic metric 
defined by $r$-dressing the cylinder with 
$h_+\in\Lambda_r^+\LieSL(2,\Bcc)_\sigma$ and the
hyperelliptic curve $\FC$ defined in Section~\ref{FCintro}.
Then, by Proposition~\ref{aonFC}, the diagonal entries of
$S=h_+Ah_+\inv$ are meromorphic on $\FC$.

For arbitrary $\zeta\in\FZ$ we write $\zeta=\hphi(\lambda)A$, where $\hphi$
is odd and meromorphic on $I^{(r)}$ with a pole only at
$\lambda=0$.
If $\hphi(\lambda)$ and $\hphi^\prime(\lambda)$ have the same
principal part at $\lambda=0$, then for $\zeta=\hphi(\lambda)A$ and
$\zeta^\prime=\hphi^\prime(\lambda)A$, $\zeta-\zeta^\prime$
is holomorphic on $I^{(r)}$. Therefore, by Lemma~\ref{higherflows}, 
$\zeta$ and $\zeta^\prime$ generate the same flow.
Thus, to show that $\Psi$ is of finite type,
it is enough to construct a trivial flow for
all but a finite dimensional space of principal parts of $\hphi$.

If $N\geq g+1$, then the function
\BEQ
\tphi_N(\nu,\mu)=\nu^{-N}\mu
\EEQ
is by~\bref{merofuncFC} 
meromorphic on $\FC$, has only a pole of order $2N-1$ at
$P_0$ and is otherwise holomorphic. Obviously,
\BEQ
\tphi_N\circ I=-\tphi_N.
\EEQ
We define
\BEQ
\phi_N=\tphi_N-\tphi_N^\ast,
\EEQ
where $(\cdot)^\ast$ was defined in Section~\ref{FCinvolutions}.
Then $\phi_N$ is a meromorphic function on $\FC$ with a pole of order
$2N-1$ at $P_0$ and $P_\infty$, which is holomorphic on $\FCstar$ and satisfies
\BEQ \label{phiNodd}
\phi_N\circ I=-\phi_N,
\EEQ
and
\BEQ \label{phiNreal}
\phi_N^\ast=-\phi_N,
\EEQ
since $(\cdot)^\ast$ is an involution which, by
Proposition~\ref{FCinvolutions}, commutes with $I$.
By~\bref{phiNodd} and Lemma~\ref{aonFC}, 
\BEQ \label{phiNarational}
\mbox{\rm $\phi_N a$ is a rational function on $\CnPE$ and an even
rational function on $\ClPE$.}
\EEQ
By~\bref{phiNreal} 
and~Proposition~\ref{aonFC}, $\phi_N a$ satisfies
\BEQ \label{phiNareal}
\mbox{\rm $(\phi_N a)^\ast=-\phi_N a$ on $\CnPE$ and on $\ClPE$.}
\EEQ
By the arguments above, Proposition~\ref{aonFC} and
Theorem~\ref{betaonFC}, the functions $\Phi_N,a,b,c,\alpha$ and
$\beta$ an all be considered being defined on ${\FC^\prime}^\ast$.
By writing $\phi_N b=\frac{\beta b}{\beta a}(\phi_N a)$,
$\phi_N c=\frac{\beta c}{\beta a}(\phi_N a)$,
using~\bref{phiNarational} and b) and d') of
Theorem~\ref{nectheorem}, we get that 
$\phi_N b$ and $\phi_N c$ are odd meromorphic functions
on $\clstar$. By writing $\phi_N b=\frac ba(\phi_N
a)$, $\phi_N c=\frac ca(\phi_N a)$, we see that $\phi_N b$ and $\phi_N c$
can be extended meromorphically to $\lambda=0$.
Furthermore, by~\bref{hahbreal}, \bref{phiNareal} and
Proposition~\ref{betaonFC} we have
\BEQ \label{phiNbreal}
(\phi_N c)^\ast=\frac{(\beta c)^\ast}{(\beta a)^\ast}(\phi_N a)^\ast=
-\frac{-\beta b}{-\beta a}(\phi_N a)=-\phi_N b.
\EEQ
In particular, $\phi_N b$ and $\phi_N c$ can be extended
meromorphically to $\lambda=\infty$, whence
\BEQ \label{phiNbrational}
\mbox{\rm $\phi_N b$ and $\phi_N c$ are odd rational functions 
on $\ClPE$.}
\EEQ
Since $a^2$, $b^2$ and $c^2$ are rational functions of $\nu$, we can
choose a polynomial $\tilde{\varphi}(\nu)$, $\tilde{\varphi}(0)\neq0$,
such that $\tilde{\varphi}^2a^2$,
$\tilde{\varphi}^2b^2$ and $\tilde{\varphi}^2c^2$ are holomorphic on
$\cnstar$.
Then $\tilde{\varphi}^\ast$ is a rational function with a pole only at
$\nu=0$. If we define 
\BEQ
\varphi=\tilde{\varphi}\tilde{\varphi}^\ast
\EEQ
then the function $\varphi(\nu)$ is rational, real on $S^1$ and
$\varphi$, $\varphi^2a^2$, $\varphi^2b^2$ and $\varphi^2c^2$ have no poles on
$\cnstar$. We set $\hphi_N=\varphi\phi_N$ on $\FCstar$.
Then, since $\phi_N$ has no
poles on $\FCstar$,~\bref{phiNarational} and~\bref{phiNbrational} yield
\BEQ \label{hphiNa}
\mbox{\rm $\hphi_N a$ is an even rational function of $\lambda$,
which is holomorphic on $\cstar$,}
\EEQ
\BEQ \label{hphiNb}
\mbox{\rm $\hphi_N b$ and $\hphi_N c$ are odd rational functions of
$\lambda$, which are holomorphic on $\cstar$.}
\EEQ
Furthermore, since $\varphi$ is real on $S^1$,~\bref{phiNareal}
and~\bref{phiNbreal} give
\BEQ \label{hphiNreal}
(\hphi_N a)^\ast=-\hphi_N a,\kern1cm\hphi_N c=-(\hphi_N
b)^\ast\kern1cm\mbox{\rm on $\cnstar$ and on $\clstar$}.
\EEQ
By~\bref{hphiNa} and~\bref{hphiNb}, the matrix function
\BEQ
(\hphi_N S)(\lambda)=\tmatrix{\hphi_N a}{\hphi_N b}{\hphi_N
c}{-\hphi_N a},
\EEQ
is rational and holomorphic on
$\clstar$ and satisfies the twisting condition~\bref{DPWtwistcond}. 
In addition, by~\bref{hphiNreal}, $(\hphi S)^\ast=-\hphi S$, i.e.,
\BEQ \label{hphiNS}
\mbox{\rm 
$\hphi_N S$ is in the Lie algebra $\Lambda_r\Liesu(2)_\sigma$ and
extends holomorphically to $\cstar$.}
\EEQ
The functions $\tilde{\varphi}$ and $\varphi$ constructed above are
independent of $N$. As in Remark~\ref{meroonFC}, we
identify $\varphi$ with an $I$-invariant meromorphic function
on $\FC$. Let $\kappa\in\Bnn$ be the degree of the
polynomial $\tilde{\varphi}(\nu)$.
Then, since $\varphi(0)\neq0$,
$\varphi$ has a pole of even order $2\kappa$, at $P_0$.
We can choose the polynomial $\tilde{\varphi}$ such that $\kappa$ is minimal.

The function $\hphi_N a$ is rational and even in $\lambda$ and has no
poles on $\cstar$.
Since, by Proposition~\ref{aonFC} and Lemma~\ref{aonFC},
$a(\lambda)\in\FA_r^+$ is odd, we get that
$\hphi_N(\lambda)=\frac{\hphi_N a}{a}$ is
an odd meromorphic function on $I^{(r)}$, which, since $\phi_N$ is
holomorphic on ${\FC^\prime}^\ast$, has only a
pole at $\lambda=0$. Therefore, by Remark~1 in
Section~\ref{higherflows}, $\zeta=\hphi_N A\in\FZ$
defines a higher flow acting on the periodic surface $\Psi$.
The pole of $\hphi_N(\lambda)$ at $\lambda=0$ is of degree $2(\kappa+N)-1$.
Since $\kappa<\infty$ is independent of $N$ and since $N$ can take
all but a finite number of integer values, it only remains to show, that
each of the $\hphi_N$ generates a trivial flow.

We have
\BEQ
h_+e^{t\hphi_NA}=h_+e^{t\hphi_N A}h_+\inv h_+=e^{t\hphi_N S}h_+.
\EEQ
By~\bref{hphiNS}, 
for all $t\in\Brr$, the matrix $U(t,\lambda)=e^{t\hphi_N S}$ is in the
twisted loop group $\Lambda_r\LieSU(2)_\sigma$ and can be extended
holomorphically to $\lambda\in\cstar$. This yields
$h_+\#(t\hphi_NA)=U_0(t)h_+$ with some unitary $\lambda$-independent
matrix $U_0(t)$, by the uniqueness of the Iwasawa decomposition.
Thus, by~\bref{Zaction} and~\bref{Finitial},
$F(z,t,\lambda)=U_0(t)F(z,\lambda)U_0(t)\inv$, and the flow
generated by $\zeta_N=\hphi_NA$ is trivial. This shows, that
$\FZ^\prime$ has at most codimension $\kappa+g<\infty$ in $\FZ$.
Thus, $\Psi$ is a finite type surface.
\QED

\myremark{}
Theorem~\ref{finitetype} can also be obtained 
from more general results for integrable systems. We want to
outline, how this works:\\
By $\FG_r$, $0<r<1$, we denote the group of maps from $C_r$ to
$\LieSL(2,\Bcc)$ which satisfy~\bref{DPWtwistcond} and can be
continued to holomorphic functions on $I^{(r)}\setminus\{0\}$.
For $0<r^\prime\leq r<1$, let
$\rho^r_{r^\prime}:\FG_r\mapsto\FG_{r^\prime}$ be the restriction of
maps. Using these homomorphisms of Lie groups, define $\FG_0$ as the direct
limit of the groups $\FG_r$ for $r\rightarrow0$.
Let $\FG_0^+\subset\FG_0$ be the subgroup of maps which extend
holomorphically to $\lambda=0$ and let $\FG_0^-$ be the subgroup of
maps which can be continued holomorphically to $\lambda=\infty$ and
take the value $I\in\LieSL(2,\Bcc)$ there. Clearly,
$\FG_0^+\cap\FG_0^-=\{I\}$.
As was proved in \cite[Section~2.2]{HaScSc:1}, the multiplication map
defines a diffeomorphism of $\FG_0^-\times\FG_0^+$ onto an open, dense subset
of $\FG_0$. By arguments similar to those in~\cite{Mc:1} one can show,
that $\FG_0$ admits a second splitting: Let $\FG_0^r\subset\FG_0$ be
the subgroup of maps which can be continued holomorphically to
$\cstar$ and take values in $\LieSU(2)$ on $S^1$. Then multiplication
$\FG_0^r\times\FG_0^+\rightarrow\FG_0$ is a diffeomorphism onto.
Since the extended frames are, by Lemma~\ref{DPWloopgroups}, all in
$\FG_0^r$ and since the dressing action depends only on the
equivalence class of $h_+$ in $\FG_0^+$, we can use the group $\FG_0$
with the splittings above, instead of
$\Lambda_r\LieSL(2,\Bcc)_\sigma$, $0<r<1$, to construct all
CMC-immersions in the $r$-dressing orbit of the cylinder. Also, the
generators $\zeta$ of the higher flows in Section~\ref{higherflows} can, by
the second statement in Lemma~\ref{higherflows}, w.l.o.g.\ be chosen such that
$e^{t\zeta}\in\FG_0^-$ for all $t\in\Brr$.
Finally, to merge both splittings into one, we can use the
classical double construction presented in~\cite[Kapitel~4.2]{Ha:1}
(compare also the analogous construction for the groups
$\Lambda_r\LieSL(2,\Bcc)_\sigma$ presented 
in~\cite{DoWu:1,Mc:1,BuPe:1,DoMcPeWu:1}). From this,
the group $\FG_0$ can be recovered by a real reduction similar to the
one used in~\cite{DoWu:1} and~\cite{BuPe:1}.
With these settings, our description of solutions of the
$\sinh$-Gordon equation fits into the framework of \cite{HaScSc:1} and
Theorem~\ref{finitetype} can be obtained from~\cite[Theorem~4.16]{HaScSc:1}.

However, for the reader's convenience, we have presented above
a direct proof of Theorem~\ref{finitetype}.

\section{Algebro-geometric description of surfaces with periodic
metric} \label{tori} \message{[tori]}
For a CMC-immersion with periodic metric, we defined in 
Section~\ref{FCintro} a nonsingular hyperelliptic curve $\FC$.
In this section, we will show, that $\FC$ allows us to express the periodicity
conditions for CMC-immersions stated in Theorem~\ref{nectheorem} and
Theorem~\ref{sufftheorem} in terms of algebro-geometric data.

If the CMC-immersion $\Psi$ under consideration does not only have a
periodic or doubly periodic metric, but even a compact image in
$\threespace$, i.e., if $\Psi$ is a CMC-torus, then
we will reproduce the classification of CMC-tori in terms of
algebro-geometric data as given in \cite{PiSt:1} and~\cite{Bo:1}.

We would like to point out that this classification refers to the
generic case. As mentioned in the introduction, a discussion of the
singular tori would be very interesting.

\newsection{} \label{cycles}
We will first reformulate the statement of Theorem~\ref{nectheorem} in
terms of algebro-geometric data. We start with the same
assumptions as in Section~\ref{periodicsurfaces}: Let
$\Psi:\Bcc\rightarrow\threespace$ be a CMC-immersion, such that
$\Sym(\Psi)$ contains a nontrivial element $q\in\cstar$.
Then we define the hyperelliptic curve $\FC$ as in Section~\ref{FCintro}.

We introduce a standard homotopy
basis for $\FC$ which is adapted to the $\FCinvolution$-symmetry of
$\FC$ stated in Proposition~\ref{FCinvolutions}.
Let $a_1,\ldots,a_g,b_1,\ldots,b_g$, $g$ the genus
of $\FC$, be a canonical basis of $H_1(\FC,\Bii)$, such that the
intersection numbers are given by
\BEQ
a_ia_j=0,\kern1cm
b_ib_j=0,\kern1cm a_ib_j=\delta_{ij},\kern1cm i,j=1,\ldots,g.
\EEQ
For the cycles $a_k$ we choose (see~\cite[VII.1.1]{FaKr:1}, where
$a_j\leftrightarrow b_j$ compared to our conventions)
\BEQ
a_k=\gamma_k-I\circ\gamma_k,
\EEQ
where $\gamma_k$ is a curve joining 
the branchpoints over $\nu_{2k-1}$ and $\nu_{2k}$, which satisfies
$\FCinvolution\circ\gamma=-\gamma$. Then
\BEQ \label{anachI}
I\circ a_k=-a_k
\EEQ
and
\BEQ \label{aktau}
\FCinvolution\circ a_k=-a_k,
\EEQ
since $\FCinvolution$ and $I$ commute.
I.e., the cycles $a_k$ are up to orientation invariant under
$\FCinvolution$ and $I$. In addition, we can choose $b_k$ such that
\BEQ \label{bktau}
\FCinvolution\circ b_k=b_k-a_k+\sum_{j=1}^ga_j.
\EEQ

\newsection{} \label{pprime}
Now we investigate the function $p\in\FA_r$ defined in~\bref{Heq}. It
is an odd meromorphic function on $I^{(r)}\subset\Bcc_\lambda$. In
general, $p$ cannot be continued to a meromorphic function on
$\FCstar$. However, the differential $\diff p=p^\prime\diff\lambda$
can. For this we note that $\alpha=\cosh(p)$ and $\beta=\sinh(p)$ on
$I^{(r)}\setminus\{0\}\subset\Bcc_\lambda$ by Theorem~\ref{nectheorem}. 
Then, since $\alpha^2-\beta^2=1$, $\diff
p=\alpha\diff\beta-\beta\diff\alpha$.

Let $\alpha$ and $\beta$ be the holomorphic functions on $\FCstar$
defined in Theorem~\ref{betaonFC}. Then we define the $1$-form
\BEQ \label{dprepr}
\mbox{\rm $\omega=\alpha\diff\beta-\beta\diff\alpha$ on $\FCstar$}.
\EEQ
Clearly,
\BEQ \label{omegahol}
\mbox{\rm $\omega$ is a holomorphic $1$-form on $\FCstar$.}
\EEQ

\mylemma{} {\em
a) In the local coordinate $\lambda=\lambda(\nu,\mu)$ around $P_0$ we
have $\omega=p^\prime(\lambda)\diff\lambda$.

b) $\rho^\ast\diff\lambda={\pi^\prime}^\ast\diff p$, where $\diff p$
is defined on $\clstar$.
}

\Proof
a) The local coordinate $\lambda(\nu,\mu)$ around $P_0\in\FC$ is the
inverse map of
\BEQ
\lambda\rightarrow(\lambda,\tmu(\lambda))
\rightarrow(\lambda^2,\lambda\tmu(\lambda)),
\EEQ
where $\tmu(\lambda)=\left(\prod_{j=1}^{2g}(\lambda-\sqrt{\nu_j})
(\lambda+\sqrt{\nu_j})\right)^\frac12$. Therefore, a) follows from
b). But $\rho^\ast\omega=\alpha(\lambda)\diff(\beta\circ\rho)
-(\beta\circ\rho)\diff\alpha(\lambda)$ and
$\beta\circ\rho=\frac{(\beta a)\circ\rho}{a\circ\rho}=\frac{(\beta
a)(\lambda)}{a(\lambda)}$, where we have used the definition of $\beta$
in Theorem~\ref{betaonFC}. For small $\lambda$, however, $(\beta
a)(\lambda)=\beta(\lambda)a(\lambda)$, whence
$\beta\circ\rho=\beta(\lambda)$ and
$\rho^\ast\omega=p^\prime(\lambda)\diff\lambda$ follows. Since in the
coordinates under consideration, $\pi^\prime(\lambda,\tmu)=\lambda$,
$(\pi^\prime)^\ast\diff p=\diff
p(\lambda)=p^\prime(\lambda)\diff\lambda$.
\QED
 
\mytheorem{} {\em 
The differential $\omega$ is a meromorphic Abelian differential of the
second kind on $\FC$ which is holomorphic
on $\FCstar$. It has poles of second order at $P_0$ and $P_\infty$.
In local coordinates around $P_0$ and $P_\infty$ it is of the form
\BEQ \label{dpsingbehaviour1}
\mbox{\rm $\omega=-q\lambda^{-2}\diff\lambda+s_+(\lambda)$, $s_+$
locally holomorphic around $P_0$,}
\EEQ
\BEQ \label{dpsingbehaviour2}
\mbox{\rm $\omega
=\qquer(\lambda\inv)^{-2}\diff(\lambda\inv)+s_-(\lambda)$, $s_-$ 
locally holomorphic around $P_\infty$}.
\EEQ
Furthermore,
\BEQ \label{dpak}
\int_{a_k}\omega=0,\kern5mm k=1,\ldots,g.
\EEQ
The differential $\omega$ is uniquely determined by these
properties. In addition, we have
\BEQ \label{dpimaginary}
\FCinvolution^\ast\omega=-\overline{\omega}.
\EEQ
}

\Proof
From~\bref{omegahol} we know that $\omega$ is a holomorphic $1$-form on
$\FCstar$. The local description of $\omega$ at $P_0$ shows that
$\omega$ has a meromorphic extension to $P_0$.

By Proposition~\ref{betaonFC} and~\bref{FCstardef},
$\alpha\circ\FCinvolution=\overline{\alpha}$ and
$\beta\circ\FCinvolution=-\overline{\beta}$. This yields
$\FCinvolution^\ast\diff\alpha=\overline{\diff\alpha}$,
$\FCinvolution^\ast\diff\beta=-\overline{\diff\beta}$.
Equation~\bref{dpimaginary} now follows from~\bref{dprepr}. Since
$\omega$ is meromorphic around $P_0$ it can by~\bref{dpimaginary} also be
continued meromorphically to $P_\infty=\FCinvolution(P_0)$. This shows,
that $\omega$ can be extended to an Abelian differential on $\FC$.

Locally around $P_0$, $\lambda$ is a local coordinate on $\FC$.
We have, by~\bref{Heq}, that
\BEQ
\omega(\lambda)=\der{p(\lambda)}{\lambda}\diff\lambda
=(-\lambda^{-2}q-\qquer+\der{f_+(\lambda)}{\lambda})\diff\lambda,
\EEQ
where $f_+$ and therefore also
$s_+=(-\qquer+\der{f_+(\lambda)}{\lambda})\diff\lambda$ 
is defined and holomorphic around $\lambda=0$.
This proves~\bref{dpsingbehaviour1}.
Locally around $P_\infty$, \bref{dpsingbehaviour1}
and~\bref{dpimaginary} yield
\BEQ
\omega=-\overline{\FCinvolution^\ast\omega}
=-\overline{\FCinvolution^\ast(-q\lambda^{-2}\diff\lambda+s_+)}
=\qquer\lambda^2\diff(\lambda\inv)+s_-
\EEQ
where $s_-=-\overline{\FCinvolution^\ast s_+}$ is locally defined and
holomorphic around $P_\infty$. This proves~\bref{dpsingbehaviour2}.
Since $\lambda\inv$ is a local coordinate around $P_\infty$ on $\FC$,
this also shows, that $\omega$ has poles of second order at
$P_\infty$ and $P_0$, and that the residue of $\omega$ vanishes at either
singularity, whence it is an Abelian differential of the
second kind. As such, it is by \cite[Prop.~III.3.3]{FaKr:1}
uniquely determined by~\bref{dpsingbehaviour1}--\bref{dpak}.

To finish the proof, it therefore suffices to verify~\bref{dpak}. To
this end we rewrite~\bref{dprepr} using $\alpha^2-\beta^2=1$ as
\BEQ \label{logderiv}
\mbox{\rm $\omega=\frac{\diff(\alpha+\beta)}{\alpha+\beta}
=\diff(\ln(\alpha+\beta))$ on $\FCstar$.}
\EEQ
Since $\alpha+\beta$ is a holomorphic function without zeroes on $\FCstar$, we
thus obtain that the integral of $\omega$ over a closed cycle on
$\FCstar$ is always an integer multiple of $2\pi i$. In particular, each
such integral is purely imaginary or zero. By~\bref{aktau}
and~\bref{dpimaginary}, we have
\BEQ
\int_{a_k}\omega=-\int_{\FCinvolution_\ast a_k}\omega
=-\int_{a_k}\FCinvolution^\ast\omega
=\int_{a_k}\overline{\omega}=\overline{\int_{a_k}\omega}.
\EEQ
Thus, the integral over $a_k$ is real. This is only possible if
the integral vanishes.
\QED

\separate
Let $\Omega_1$ be an Abelian differential of the second kind on
$\FC$ which is in local coordinates around $P_\infty$ given by
\BEQ \label{dOsingbehaviour}
\Omega_1=-(\lambda\inv)^{-2}\diff(\lambda\inv)+\mbox{\rm a
locally holomorphic differential around $P_\infty$}
\EEQ
and is holomorphic on $\FC\setminus\{P_\infty\}$. Such a differential
exists by~\cite[Theorem~II.5.1]{FaKr:1}. It is unique up to addition
of a holomorphic differential on $\FC$. Thus, 
see \cite[Corollary~III.3.4]{FaKr:1}, if we normalize $\Omega_1$ by
\BEQ \label{dO1ak}
\int_{a_k}\Omega_1=0,\kern5mm k=1,\ldots,g,
\EEQ
then $\Omega_1$ is determined uniquely by~\bref{dOsingbehaviour}
and~\bref{dO1ak}. By~\bref{anachI}, the differential
\BEQ \label{O2def}
\Omega_2=\overline{\FCinvolution^\ast\Omega_1}
\EEQ
has also vanishing $a$-periods and
a single pole at $P_0$. It is given in local coordinates around $P_0$ by
\BEQ
\Omega_2=-\lambda^{-2}\diff\lambda+\mbox{\rm a locally
holomorphic differential around $P_0$}.
\EEQ
Since Abelian differentials are uniquely determined by their
$a$-periods and the principal parts
at the singularities~\cite[Prop.~III.3.3]{FaKr:1}, we have
\BEQ \label{pOmega}
\omega=-\qquer\Omega_1+q\Omega_2
=-\qquer\Omega_1+q\overline{\FCinvolution^\ast\Omega_1}.
\EEQ
In addition $I^\ast\Omega_1$ and $-\Omega_1$ have the same
pole divisors and both
satisfy~\bref{dO1ak} and~\bref{dOsingbehaviour}. Therefore,
\BEQ \label{dOI}
I^\ast\Omega_1=-\Omega_1.
\EEQ
Let us define
\BEQ \label{UVdef}
U_k=\int_{b_k}\Omega_1,\kern3cm V_k=\int_{b_k}\Omega_2.
\EEQ
From~\bref{bktau} it follows, that
\BEQ \label{VkUkquer}
V_k=\int_{b_k}\Omega_2
=\int_{b_k}\overline{\FCinvolution^\ast\Omega_1}
=\int_{\FCinvolution_\ast b_k}\overline{\Omega_1}
=\overline{\int_{b_k}\Omega_1}=\overline{U_k}.
\EEQ
By collecting the results above, we get:

\myprop{} {\em
The integrals $U_k$ of $\Omega_1$ over the $b$-cycles satisfy
\BEQ \label{UVcond}
\Imag(q\overline{U_k})=\pi m_k,\kern5mm m_k\in\Bii,\kern1cm k=1,\ldots g.
\EEQ
}

\Proof
By~\bref{pOmega} and~\bref{VkUkquer}, we have
\BEQ \label{bcycleUV}
\int_{b_k}\omega=-\qquer U_k+q\overline{U_k}=2i\Imag(q\overline{U_k}).
\EEQ
As in the proof of Lemma~\ref{pprime}, we use the fact that the integral
of $\omega$ over a closed cycle is an integer multiple of $2\pi
i$. This together with~\bref{bcycleUV} proves~\bref{UVcond}.
\QED

\newsection{} \label{necclosing}
Next we investigate the closing conditions in
Theorem~\ref{spaceclosing}.
Since $\beta^2$ is an even function in $\lambda$, we can write it as a
function of $\nu=\lambda^2$. By Lemma~\ref{FCintro}, $a^2(\nu)$ has no
zero of odd order on $S^1$. Therefore, by d') in
Theorem~\ref{nectheorem}, also $\beta^2=\frac{(\beta a)^2}{a^2}$ has 
no zeroes of odd order on $S^1$. As a function
of $\lambda$, $\beta^2$ has a zero of order $2n$, $n>0$, at
$\lambda_0\in S^1$ iff as a function of $\nu$ it has a zero of
order $2n$ at $\nu_0=\lambda_0^2$.
Since $\nu$ is a local coordinate around every point on $\FS$, we have that
$\beta^2$ has a zero of order $2n$ at $\nu_0$ iff $\beta$ has a zero
of order $n$ at each of the two covering points of $\nu_0$ on
$\FC$. We have proved the following

\mylemma{} {\em
Let $\lambda_0\in S^1$ and let $P_1$ and $P_2$ be the two points on
$\FC$, such that $\nu(P_1)=\nu(P_2)=\lambda_0^2$. Then we have
\begin{description}
\item[a)] $\beta^2(\lambda)$ vanishes at $\lambda_0$ iff $\beta$ vanishes at
$P_0$ and $P_1$.
\item[b)] $\beta^2(\lambda)$ vanishes at least to fourth order at
$\lambda_0$ iff $\beta$ vanishes at least to second order at $P_0$
and $P_1$.
\end{description}
}

\separate
Since $\beta(P_0)=-\beta(P_1)$, we have that in
both statements of Lemma~\ref{necclosing} we can replace ``and'' by ``or''.
We can restate Theorem~\ref{spaceclosing} in terms of $\omega$:

\mytheorem{} {\em
Let $\Psi:\Bcc\rightarrow\threespace$ be a CMC-immersion in the
$r$-dressing orbit of the cylinder with associated family
$\{\Psi_\lambda,\lambda\in S^1\}$. Let
$q\in\Sym(\Psi)=\Sym(\Psi_\lambda)$, $q\neq0$, and let $\chi$ be
defined as in Theorem~\ref{nectheorem}. Let $\FC$ be the
hyperelliptic curve defined in Section~\ref{FCintro} and let $\omega$
be the Abelian differential on $\FC$ which is uniquely defined by
Theorem~\ref{pprime}. For $\lambda_0\in S^1$ denote by $P_1(\lambda_0)$
and $P_2(\lambda_0)$ the covering points of $\nu_0=\lambda_0^2$ on $\FC$.
Then the following are equivalent:
\begin{description}
\item[1.] $q\in\Per(\Psi_{\lambda_0})$.
\item[2.] $\beta$ vanishes at least to second order in
$P_1(\lambda_0)$ and $P_2(\lambda_0)$.
\item[3.] The form $\omega$ has a zero at $P_1(\lambda_0)$ or
$P_2(\lambda_0)$ and there
exists a curve $\gamma$ on $\FCstar$ connecting $P_1(\lambda_0)$ and
$P_2(\lambda_0)$ such that
\BEQ \label{secondclosing}
\int_\gamma\omega=2\pi im,\kern5mm m\in\Bii.
\EEQ
\end{description}
}

\Proof
1.$\Leftrightarrow$2. follows from Theorem~\ref{spaceclosing} and
Lemma~\ref{necclosing}.

2.$\Rightarrow$3. Using~\bref{logderiv} we get 
for each curve on $\FCstar$ connecting $P_1$ and $P_2=I(P_1)$:
\BEQ \label{intomega}
\int_\gamma\omega=\ln(\alpha+\beta)|_{P_1}^{P_2}
=\ln((\alpha+\beta)(P_2))-\ln((\alpha+\beta)(P_1))
\EEQ
where we have continued $\ln(\alpha+\beta)$ analytically along
$\gamma$. Since $\beta(P_1)=\beta(P_2)=0$,
\BEQ
\int_\gamma\omega=\ln(\alpha(P_2))-\ln(\alpha(P_1)),
\EEQ
which by~\bref{463} is an integer multiple of $2\pi i$.
Thus,~\bref{secondclosing} holds. In addition, if $\beta$ vanishes
to second order at $P_1$ and $P_2$, then also $\diff\beta$ vanishes at
both points. Therefore, $\omega=\alpha\diff\beta-\beta\diff\alpha$ has
a zero at $P_1$ and $P_2$.

3.$\Rightarrow$2. Let $\gamma$ be a curve on $\FCstar$ connecting
$P_1$ and $P_2$ for which~\bref{secondclosing} holds.
Then by~\bref{intomega} we get
\BEQ
(\alpha-\beta)(P_1)=(\alpha+\beta)(P_1)
\EEQ
from which with Proposition~\ref{betaonFC} $\beta(P_1)=-\beta(P_2)=0$
follows. If in addition
$\omega(P_1)=0$, then since $\beta(P_1)=0$, we get $\alpha(P_1)=1$ and
$\diff\beta(P_1)=0$. Thus $\beta$ vanishes to second order at $P_1$
and therefore also at $P_2$.
\QED

\myremark{}
1. It follows immediately from~\bref{Fdressingdef}, that
\BEQ
F(z,-\lambda)=F(-z,\lambda),
\EEQ
i.e., $\Psi_\lambda$ and $\Psi_{-\lambda}$ differ only by a
coordinate transformation. This explains, why the
conditions 2.~and 3.~in Theorem~\ref{necclosing} are invariant under
the substitution $\lambda_0\rightarrow-\lambda_0$.

2. The proof of Theorem~\ref{necclosing} shows, 
that in the third statement of the theorem, ``there
exists a curve $\gamma$ on $\FCstar$'' can be replaced by 
``for all curves $\gamma$ on $\FCstar$''.

\newsection{} \label{Bobenkotori}
If $\Psi$ is a periodic CMC-immersion with associated family
$\{\Psi_\lambda;\lambda\in S^1\}$.
then by Remark~\ref{DPWperiodic}, $q\in\Sym(\Psi_\lambda)$
for all $\lambda\in S^1$. In Sections~\ref{FCintro}, \ref{cycles} 
and~\ref{pprime} we have
introduced a nonsingular hyperelliptic curve $\FC$, a canonical
homotopy basis $a_1,\ldots,a_g,b_1,\ldots,b_g\in H_1(\FC,\Bii)$ and an
Abelian differential $\Omega_1$ of the second kind on
$\FC$. If $q\in\Sym(\Psi)$, then these data satisfy 
Proposition~\ref{pprime}.
If for some $\lambda_0\in S^1$, $q\in\Per(\Psi_{\lambda_0})$, then in
addition the statement of Theorem~\ref{necclosing} holds.
In the special case that $\Psi_{\lambda_0}$ is a CMC-torus, we get

\mytheorem{} {\em
Let $\Psi:\Bcc\rightarrow\threespace$ be a periodic CMC-immersion with
associated family $\{\Psi_\lambda;\lambda\in S^1\}$. We define the
nonsingular hyperelliptic curve $\FC$ as in Section~\ref{FCintro}.
Assume, that
$\Psi_{\lambda_0}$ is a CMC-torus. Let $q_j=\frac12Y_j+i\frac12X_j$, 
$j=1,2$, be the
generators of $\Per(\Psi_{\lambda_0})$. We introduce a homotopy basis
$a_1,\ldots,a_g,b_1,\ldots b_g$, the Abelian differential
$\Omega_1$ and $U_k=\alpha_k+i\beta_k$ as in
Sections~\ref{cycles} and~\ref{pprime}.
Then there exists a curve $\gamma$
on $\FCstar$
connecting the points $P_1(\lambda_0)$ and $P_2(\lambda_0)$ for which
$\nu(P_i(\lambda_0))=\lambda_0^2$, such that for 
$X_1,Y_1,X_2,Y_2$, $\alpha_k,\beta_k$,
$k=1,\ldots,g$, defined above and $c_1,c_2$ defined by
\BEQ
c_1+ic_2=\frac12\int_\gamma\Omega_1
\EEQ
\begin{description}
\item[1.] the matrix
\BEQ \label{intmatrix}
\frac{1}{2\pi}\tmatrix{X_1}{Y_1}{X_2}{Y_2}
\tmatrix{2c_1}{\alpha_1,\ldots,\alpha_g}{-2c_2}{-\beta_1,\ldots,-\beta_g}
\EEQ
has integer entries and
\item[2.] $\Omega_1$ vanishes at $P_1(\lambda_0)$ 
(and therefore also at $P_2(\lambda_0)$).
\end{description}
}

\Proof
Let $\Psi_{\lambda_0}$ be a CMC-torus, such that $q_1$ and $q_2$
generate $\Per(\Psi_{\lambda_0})$. If we define $\alpha_k,\beta_k$,
$k=1,\ldots,g$, as above then by Proposition~\ref{pprime} we have for
$j=1,2$:
\BEQ
\pi m_k=\Imag(q_j\overline{U_k})
=\frac12(X_j\alpha_k-Y_j\beta_k),
\EEQ
where $m_k\in\Bii$. This shows, that the last $g$ columns
in the matrix product of~\bref{intmatrix} have integer entries.

In addition, if we set $V=\int_\gamma\Omega_1$, then by
Theorem~\ref{necclosing}, we get for $j=1,2$:
\BEQ
2i\Imag(q_j\overline{V})-\overline{q_j}V+q_j\overline{V}
=2i(X_jc_1-Y_jc_2)=2\pi im_j,\kern1cm m_j\in\Bii.
\EEQ
This shows, that the first column of~\bref{intmatrix} has integer
entries, proving the first part of the statement.

Finally, by Theorem~\ref{necclosing} and~\bref{pOmega},
we get that $-\overline{q_1}\Omega_1+q_1\overline{\Omega_1}$ 
and $-\overline{q_2}\Omega_1+q_2\overline{\Omega_1}$ vanish
at $P_1(\lambda_0)$ or $P_2(\lambda_0)$, and therefore, by~\bref{dOI},
at both points. Since
$q_1$ and $q_2$ are linearly independent as vectors in $\twospace$,
it follows that $\Omega_1$ vanishes at $P_1(\lambda_0)$.
\QED

\myremark{}
This reproduces the conditions~(6.13) and (6.14) of \cite{Bo:1} for
CMC-tori.

\newsection{}  \label{FCconstruction}
In the last two sections we have reformulated the necessary conditions
of Theorem~\ref{nectheorem} for periodic surfaces in terms of
algebro-geometric data on $\FC$. For CMC-tori, the result is the same as in
\cite[Theorem 6.1]{Bo:1}, as expected.
It remains to construct periodic surfaces from algebro-geometric data.

We start with the following data:
\begin{description}
\item[1.] Let $g\in\Bnn$ and let $\nu_{2k-1}\in\Bcc$, $k=1,\ldots,g$,
be mutually distinct, such that $0<|\nu_{2k-1}|<1$,
\item[2.] let $q\in\cstar$.
\end{description}
Define the nonsingular hyperelliptic curve $\FC$ by
\BEQ \label{suffFCdef}
\mu^2=\nu\prod_{k=1}^{2g}(\nu-\nu_k),\kern1cm
\nu_{2k}=\overline{\nu_{2k-1}}\inv,\kern5mm k=1,\ldots,g.
\EEQ
Then $g$ is the genus of $\FC$ and $\FC$ has no branchpoints over the
unit circle. For this surface, Theorem~\ref{FCintro} holds.
The surface $\FC$ admits the
involution $\FCinvolution$ defined in~\bref{tauonFC}. As in
Section~\ref{FCinvolutions}, we define for a scalar function $f$ on $\FC$,
$f^\ast=\overline{f\circ\FCinvolution}$.
We denote by $P_0$ and
$P_\infty$ the branchpoints of $\FC$ over $\nu=0$ and $\nu=\infty$.
We also define a standard homotopy basis
$a_1,\ldots,a_g,b_1,\ldots,b_g$, such that~\bref{aktau}
and~\bref{bktau} hold.
Let $\Omega_1$ be the uniquely defined Abelian differential on $\FC$ which
is holomorphic except at $P_\infty$ and satisfies~\bref{dO1ak}
and~\bref{dOsingbehaviour}.
We also define $\Omega_2$ by~\bref{O2def}.
We denote the integral of $\Omega_1$ over the cycle $b_k$ by $U_k$.

We now assume, that with these definitions the equation
\BEQ \label{qsuffcond}
\Imag(q\overline{U_k})=\pi m_k,\kern5mm m_k\in\Bii,\kern1cm k=1,\ldots g,
\EEQ
is satisfied.
{}From these data we will construct a family of CMC-immersions
which are periodic w.r.t.\ translation by $q$.

For later use we also introduce the hyperelliptic curve $\FC^\prime$
defined by
\BEQ
\tmu^2=\prod_{k=1}^{2g}(\lambda-\sqrt{\nu_k})(\lambda+\sqrt{\nu_k}).
\EEQ
The discussion in Sections~\ref{meroonFC}--\ref{FCprime}
applies to $\FC$ and $\FC^\prime$.
In particular, we have the local coordinate $\lambda=\lambda(\nu,\mu)$
at $P_0$ on $\FC$ (see also the proof of Lemma~\ref{pprime}).

\newsection{} \label{alphabetaconstruction}
Next we define the Abelian differential $\omega$ by~\bref{pOmega}.
By the definition of $\Omega_1$ and $\Omega_2$, it
satisfies~\bref{dpimaginary} and~\bref{dpak}. It has poles at $P_0$
and $P_\infty$. Its behaviour at the
singularities is given by~\bref{dpsingbehaviour1},\bref{dpsingbehaviour2}.
Let $I$ be the hyperelliptic involution on $\FC$.
Since $-I^\ast\omega$ also
satisfies~\bref{dpsingbehaviour1}--\bref{dpak}, we have
\BEQ \label{dpI}
I^\ast\omega=-\omega.
\EEQ
We define the multi-valued holomorphic function $p$ on
$\FCstar=\FC\setminus\{P_0,P_\infty\}$ by
\BEQ
p(\nu,\mu)=\int_{\gamma(\nu,\mu)}\omega,
\EEQ
where $\gamma(\nu,\mu)$ is a curve connecting the branchpoint over
$\nu_1$ with $(\nu,\mu)$. By~\bref{dpak}, \bref{qsuffcond}
and \bref{bcycleUV}, $p$ is
at every point of $\FCstar$ defined up to an integer multiple of $2\pi i$.
Moreover, every branch of $p$ is locally meromorphic at each point of
$\FCstar$. The functions
\BEQ
\alpha=\cosh(p),\kern3cm\beta=\sinh(p)
\EEQ
are holomorphic and single-valued on $\FCstar$.
Clearly, $\alpha^2-\beta^2=1$ on $\FCstar$.
By~\bref{dpI}, we have modulo $2\pi i\Bii$:
\BEQ
(p\circ I)(\nu,\mu)=\int_{\gamma(\nu,-\mu)}\omega
=\int_{\gamma(\nu,\mu)}(I^\ast\omega)=-p(\nu,\mu),
\EEQ
where we have chosen $\gamma(\nu,\mu)=I\circ\gamma(\nu,-\mu)$.
By~\bref{dpimaginary}, we have modulo $2\pi i\Bii$:
\BEQ \label{podd}
(p\circ\FCinvolution)(\nu,\mu)=\int_{\gamma(\FCinvolution(\nu,\mu))}\omega=
\int_{\gamma(\nu,\mu)}(\FCinvolution^\ast\omega)
=-\int_{\gamma(\nu,\mu)}\overline{\omega}=-\overline{p(\nu,\mu)},
\EEQ
where we have chosen
$\gamma(\nu,\mu)=\FCinvolution\circ\gamma(\FCinvolution(\nu,\mu))$.
This shows, that $p^\ast=-p$ modulo $2\pi i\Bii$. Therefore, we get
\BEQ \label{alphabetatau}
\alpha^\ast=\alpha,\kern3cm\beta^\ast=-\beta
\EEQ
and
\BEQ \label{alphabetaI}
\alpha\circ I=\alpha,\kern3cm\beta\circ I=-\beta.
\EEQ
Thus, $\alpha$ and $\beta^2$ project down to holomorphic functions
$\alpha(\nu)$, $\beta^2(\nu)$ on $\cnstar$ which are real on $S^1$.
Pulling back $\omega,p,\alpha$ and $\beta$ from $\FCstar$ to
${\FC^\prime}^\ast$ and from $\cnstar$ to $\clstar$,
respectively, yields holomorphic functions $\alpha(\lambda)$ and
$\beta^2(\lambda)$, which are real on $S^1$ in view
of~\bref{alphabetatau}. Moreover, \bref{alphabetaI} and 
\bref{alphabetatau} show, that $\alpha$ and $\beta^2$ are even on
$\clstar$ and that $\beta$ is non-positive over $S^1$ by
Lemma~\ref{FCinvolutions}.

The function $p$ on ${\FC^\prime}^\ast$ can be written locally near
one of the points on $\FC^\prime$ above $\lambda=0$ in the form
$p(\lambda)=\lambda\inv q-\lambda\qquer+f_+$, where $f_+$ is 
locally holomorphic at $\lambda=0$ and, by~\bref{podd}, an odd
function of $\lambda$.

We collect these results in the following

\mylemma{} {\em
With the notation of Section~\ref{FCconstruction} we have:
Let $0<r_0<1$ be such that
$\FC^\prime$ has no branchpoints over $I^{(r_0)}$. Define
\BEQ
p(\nu,\mu)=\int_{\gamma(\nu,\mu)}
(-\qquer\Omega_1+q\overline{\FCinvolution^\ast\Omega_1}),
\EEQ
where $\gamma(\nu,\mu)$ is an arbitrary curve joining $(\nu_1,0)$ and
$(\nu,\mu)$ on $\FC$. Then $p$ is a multivalued function on $\FC$
unbranched over $\nu\in I^{(r_0)}$. Every branch of $p$ can be identified
with an odd meromorphic function on $\lambda\in I^{(r_0)}$, which is of
the form
\BEQ
p(\lambda)=\lambda\inv q-\lambda\qquer+f_+(\lambda),
\EEQ
where $f_+(\lambda)$ is odd in $\lambda$ and
holomorphic on $I^{(r_0)}$, i.e., $f_+\in\FA_{r_0}^+$.
If we define $\alpha(\lambda)=\cosh(p(\lambda))$,
$\beta(\lambda)=\sinh(p(\lambda))$, then $\alpha$ and $\beta$ satisfy
the conditions a')--c') of Theorem~\ref{sufftheorem}.
Furthermore, $\beta$ vanishes at all branchpoints of $\FC$,
\BEQ \label{betavanishes}
\beta(\nu_k,0)=0,\kern1cm k=1,\ldots,2g.
\EEQ
}

\Proof
It only remains to prove~\bref{betavanishes}. But this follows
immediately from~\bref{alphabetaI} since the branchpoints are
fixed points under the hyperelliptic involution $I$ on $\FC$.
\QED

\newsection{} \label{aconstruction}
It remains to construct rational
functions $a^2(\lambda)$, $b^2(\lambda)$, $c^2(\lambda)$,
such that the conditions a)--e) and d') of
Theorem~\ref{sufftheorem} are met. We will first construct $a^2(\lambda)$.

If the genus $g$ of $\FC$ is odd, we define
\BEQ \label{evena2}
\ha^2(\nu)=\left(\prod_{k=1}^g\nu_{2k-1}\right)
\nu^g\prod_{k=1}^{2g}(\nu-\nu_k)\inv.
\EEQ
If $g$ is even, we choose an arbitrary $\nu_0\in S^1$ and define
\BEQ \label{odda2}
\ha^2(\nu)=\nu_0\left(\prod_{k=1}^g\nu_{2k-1}\right)
\nu^{g-1}(\nu-\nu_0)^2\prod_{k=1}^{2g}(\nu-\nu_k)\inv.
\EEQ
We note, that in both cases $\ha^2$ is a rational function with a zero
of odd order at $\nu=0$ and simple poles precisely at the points
$\nu_k$, $k=1,\ldots,2g$. A simple calculation gives
\BEQ \label{tareal}
(\ha^2)^\ast(\nu)=\overline{\ha^2(\nuquer\inv)}=\ha^2(\nu).
\EEQ
In addition, $\ha^2$ has no zeroes of odd order different from
$\nu=0$. In particular, since
$\ha^2$ is real on $S^1$ with at most one zero on $S^1$, 
it is either non-negative or non-positive on
$S^1$. Therefore, there exists a real constant $0<A<\infty$, such that
$0\leq\ha^2(\nu)\leq A$ or $-A\leq\ha^2(\nu)\leq0$ for all $\nu\in
S^1$. Let us define the rational function
\BEQ \label{aa0}
a_0^2=\epsilon\ha^2,
\EEQ
where $\epsilon=\pm1$ is chosen such that $a_0^2$ is non-negative on
$S^1$.

\myprop{} {\em
Let $\FC$, $\FC^\prime$, $\beta$ and $0<r_0<1$ be defined as in
Lemma~\ref{alphabetaconstruction}. Then $a_0^2(\nu)$ constructed above
is a rational even function of $\lambda$, which is real and
non-negative on $S^1$. Furthermore,
$a_0^2(\lambda)$ restricts to the square of an odd holomorphic 
function $a_0(\lambda)\in\FA_{r_0}^+$,
such that $(\beta a_0)(\lambda)$ can be extended holomorphically to an
even function on $\clstar$.
}

\Proof
We have already shown, that $a_0^2$ is rational in $\nu=\lambda^2$ and
therefore an even rational function of $\lambda$. From~\bref{tareal}
and~\bref{aa0}, it follows, that $a_0^2$ is real
and non-negative on $S^1$. Since $a_0^2$ has
only a zero of odd order at $\nu=0$
and only poles of odd order at $\nu_k$,
$k=1,\ldots,2g$, its square root $a_0$ is a meromorphic function on
$\FC$, which satisfies $a_0\circ I=-a_0$. 
It has simple poles precisely at the branchpoints of $\FCstar$
and has zeroes of odd order at $P_0$ and $P_\infty$.
Applying the usual identifications between functions on $I^{(r_0)}$,
functions on a patch of $\FC^\prime$ above $I^{(r_0)}$ and functions
in a neighbourhood of $P_0$ on $\FC$, and also using Lemma~\ref{aonFC}
we conclude that $a_0$ is
an odd holomorphic function of $\lambda$ on $I^{(r_0)}$, i.e.,
$a_0\in\FA_{r_0}^+$. In addition, since $\beta$
is holomorphic on $\FCstar$ and, by~\bref{betavanishes},
vanishes at the branchpoints of $\FCstar$, the function
$\beta a_0$ is holomorphic on $\FCstar$ and, by~\bref{alphabetaI}
invariant under $I$. Therefore, $\beta a_0$ projects down to an even
holomorphic function $(\beta a_0)(\lambda)$ on $\cstar$.
\QED

\separate
In the following, for $x\in\Brr$ we denote by $[x]$ the greatest
integer less than or equal to $x$.  

\mylemma{} {\em
Let $\FC$, $\FC^\prime$, $\beta$ and $0<r_0<1$ be defined as in
Theorem~\ref{alphabetaconstruction}.
Let $\tf(\nu)$ be a rational function which is holomorphic on
$\cstar$, real on $S^1$,
and has a pole of order at most $[\frac{g-1}{2}]$ at $\nu=0$.
Let us define $a_0^2$ as above and set $a^2=\tf^2 a_0^2$. For each such
$\tf$, $a^2$ is an even
rational function of $\lambda$ which is real and non-negative on $S^1$.
Furthermore, the restriction of $a^2$ to $C_{r_0}$, the circle with
radius $r_0$ with center $\lambda=0$, is the
square of an odd function $a\in\FA_{r_0}^+$. In addition,
$\beta a$ extends to an even 
holomorphic function $(\beta a)(\lambda)$ on $\cstar$. 
}

\Proof
By Proposition~\ref{aconstruction}, the 
functions $\tf^2$, $a_0^2$ and therefore also $a^2$ are even
rational functions of $\lambda$, which are real and non-negative on $S^1$.
In addition, the restriction of $a_0^2$ to $C_{r_0}$
is the square of an odd function $a_0(\lambda)\in\FA_{r_0}^+$. Thus,
since $\tf$ is even in $\lambda$, and since $\tf^2a^2$ has no poles on
$I^{(r_0)}$, the restriction of $a^2$ to $C_{r_0}$ is the square of the odd
function $a(\lambda)=\tf(\lambda)a_0(\lambda)\in\FA_{r_0}^+$.
Furthermore, $\tf$, $\beta a_0$ and therefore also $\beta a=\tf\cdot(\beta
a_0)$ are even holomorphic functions of $\lambda$ on $\cstar$.
\QED

\mytheorem{} {\em
Let $\FC$, $\FC^\prime$, $\beta$ and $0<r_0<1$ be defined as in
Theorem~\ref{alphabetaconstruction}. Then there exists at least a
real $g$-parameter family of even rational functions $a^2(\lambda)$, such
that
\begin{description}
\item[a)] $a^2$ is real on $S^1$ and $0\leq a^2(\nu)<1$ for $\nu\in S^1$,
\item[b)] The restriction of $a^2$ to $C_{r_0}$ is the square of an
odd function $a(\lambda)\in\FA_{r_0}^+$,
\item[c)] $(\beta a)(\lambda)$ can be extended to an even holomorphic
function on $\cstar$.
\end{description}
}

\Proof
By Lemma~\ref{aconstruction}, to each
rational function $\tf(\nu)$ which is real on $S^1$, holomorphic on
$\cstar$, and has a pole of order at most $[\frac{g-1}{2}]$ we can construct an
even rational function $a^2(\lambda)$ which satisfies b) and c), is real
and non-negative on $S^1$ and has no poles on $S^1$. Thus, we can
always normalize $a^2$ such that $a^2(\nu)<1$ for all $\nu\in
S^1$. Then also a) is satisfied.

Since $\tf$ is real on $S^1$, we have, by Section~\ref{involutions},
$\tf^\ast=\tf$. Therefore, $\tf$ has the same pole order at $\nu=0$
and $\nu=\infty$. As a consequence, with $m=[\frac{g-1}{2}]$, $\tf$ is
of the form $\lambda^{-m}(r_0+r_1\lambda+\ldots
+r_{m-1}\lambda^{m-1})+t_0+t_1\lambda+\ldots+t_m\lambda^m$.
This implies that $\tf$ is defined by
$2[\frac{g-1}{2}]+1$ complex parameters. After imposing the reality
condition there are still
$2[\frac{g-1}{2}]+1$ real parameters left free, which determine the chosen
admissible function $\tf$ completely.
If $g$ is odd, then $2[\frac{g-1}{2}]=g-1$.
If $g$ is even, then $2[\frac{g-1}{2}]=g-2$. In the case that $g$ is
even we have the additional freedom to choose a zero $\nu_0\in S^1$ of $a_0^2$.
Altogether, in both cases this gives a real $g$-parameter family of
functions $a^2$, which satisfy a)--c).
\QED

\newsection{} \label{bcconstruction}
It remains to construct the even rational functions $b^2$ and $c^2$,
such that all conditions of Theorem~\ref{sufftheorem} are satisfied.

\mylemma{} {\em
Let $a^2(\lambda)$ and $0<r_0<1$ be as in 
Theorem~\ref{aconstruction}. Then the equation
\BEQ \label{b2eq}
b^2(b^2)^\ast=(1-a^2)^2,
\EEQ
has a solution $b^2(\lambda)$ which is rational, even in $\lambda$,
and its restriction to $C_{r_0}$ is the square of an even
function $b(\lambda)\in\FA_{r_0}^+$. 
Furthermore, we can choose $b^2$ such that $b^2$ and $(b^2)^\ast$ both
have either a simple zero
or a simple pole at each of the branchpoints $\nu_1,\ldots,\nu_{2g}$ of $\FC$
and no other poles or zeroes of odd order.
}

\Proof
With $a^2$ also $1-a^2$ is rational, even in $\lambda$, and defined and
real on $S^1$. Furthermore, $1-a^2$ has
no zeroes on $S^1$. By reordering the branchpoints
$\nu_1,\ldots,\nu_{2g}$ if necessary we can assume that $a^2$ has
simple poles at $\nu_1,\ldots,\nu_{2K}$ and zeroes of odd order at
$\nu_{2K+1},\ldots,\nu_{2g}$, where $g-[\frac{g-1}{2}]\leq K\leq g$.
Here, the lower bound for the number of poles of $a^2$ follows from
the maximum number $2[\frac{g-1}{2}]$ of zeroes of the rational
function $\tf$ in Lemma~\ref{aconstruction}.
We write
\BEQ
1-a^2=\gamma
\frac{\prod_{j=1}^{n_1}(\nu-\nu_j^{(1)})(\nu-\overline{\nu_j^{(1)}}\inv)}{
\prod_{k=1}^K(\nu-\nu_{2k-1})(\nu-\nu_{2k})},\kern1cm\gamma\in\Brr,
\EEQ
where $\{\nu_j^{(1)},\overline{\nu_j^{(1)}}\inv\}$ 
are the zeroes of $1-a^2$.
Since $(1-a^2)(\lambda=0)=(1-a^2)^\ast(\lambda=\infty)=1$, we have
$n_1=K$ and $\nu_j^{(1)}\neq0$ for all $j=1,\ldots,n_1$.
Define
\BEQ \label{b2def}
\tb^2=\prod_{k=1}^K\frac{(\nu-\nu_k^{(1)})^2}{
(\nu-\nu_{2k-1})(\nu-\nu_{2k})}
\prod_{j=K+1}^g\frac{\nu-\nu_{2j-1}}{\nu-\nu_{2j}},
\EEQ
then a straightforward computation shows that there exists a positive
real constant $\delta$, such that 
$(1-a^2)^2=\delta(\tb^2)(\tb^2)^\ast$. If we define
$b^2=\sqrt{\delta}\tb^2$ then $b^2$ is rational and
satisfies~\bref{b2eq}. Also, $b^2$ has a simple pole or a simple zero
at each of the points $\nu_1,\ldots,\nu_{2g}$
and no other poles or zeroes of odd
order. The function $(b^2)^\ast$ is up to a complex constant given by
\BEQ
\prod_{k=1}^K\frac{(\nu-\overline{\nu_k^{(1)}}\inv)^2}{
(\nu-\nu_{2k-1})(\nu-\nu_{2k})}
\prod_{j=K+1}^g\frac{\nu-\nu_{2j}}{\nu-\nu_{2j-1}},
\EEQ
and has therefore also simple poles or simple zeroes at
$\nu_1,\ldots,\nu_{2g}$ and no other poles or zeroes of odd order.

Since $b^2$ has no poles and no zeroes of odd order in
$I^{(r_0)}$, its restriction to $C_{r_0}$ is the square of a
function $b(\lambda)\in\FA_{r_0}^+$.
Since we have written $b^2$ in terms of $\nu$, it is clear that $b^2$
is an even function on $\Bcc_\lambda$. Thus, $b(\lambda)$ is even or odd. Since
$b^2$ does not vanish at $\lambda=0$, we get that $b(\lambda=0)\neq0$ whence
$b(\lambda)$ is even.
\QED

\myprop{} {\em
Let $\FC$ and $q$ be defined as in Section~\ref{FCconstruction}.
Define functions $f_+\in\FA_{r_0}^+$,
$0<r_0<1$, $p=\lambda\inv q-\lambda\qquer+f_+$, $\alpha=\cosh(p)$,
$\beta=\sinh(p)$ as in Lemma~\ref{alphabetaconstruction}.
Furthermore, let $a^2(\lambda)$ be defined as in
Theorem~\ref{aconstruction}.
Then there exist two even rational functions
$b^2(\lambda)$, $c^2(\lambda)$, such that all the conditions of
Theorem~\ref{sufftheorem} are satisfied for some $r^\prime=r=r_0$.
}

\Proof
By Lemma~\ref{alphabetaconstruction}, there exists $0<r_0<1$, such that
the functions $f_+$, $p$, $\alpha$ and $\beta$ satisfy Conditions a')--c') of 
Theorem~\ref{sufftheorem} for $r^\prime=r_0$, where the maximal radius
$r_0$ is determined by the surface $\FC$.
By a),b) in Theorem~\ref{aconstruction}, each function $a^2$
constructed in Theorem~\ref{aconstruction} satisfies Conditions a) and c) of
Theorem~\ref{sufftheorem} for $r=r_0$. Furthermore, by c) in
Theorem~\ref{aconstruction}, $(\beta a)(\lambda)$
can be extended holomorphically to $\cstar$.
We define the even rational function
$b^2(\lambda)$ as in~\bref{b2def} and use Condition b)
of Theorem~\ref{sufftheorem} to define $c^2$.
Clearly, $c^2(\lambda)$ is an even rational function which has the
same pole divisor as $b^2$ and $a^2$ and only zeroes of even
order.
Since $c^2$ has no pole and no zero of odd order on $I^{(r_0)}$, its
restriction to $C_{r_0}$ is a function $c(\lambda)\in\FA_{r_0}^+$.
This function is even or odd in $\lambda$.
Since $b^2$ does not vanish at $\lambda=\infty$, $c^2(\lambda=0)\neq0$
and $c$ is an even function in $\lambda$. This shows, that Condition
c) and d) of Theorem~\ref{sufftheorem} are satisfied.
By~\bref{b2eq}, $b^2c^2=(1-a^2)^2$. We choose the sign of the square
roots $b$ and $c$ of $b^2$ and $c^2$, such that $bc=1-a^2$. Then also
Condition e) of Theorem~\ref{sufftheorem} is satisfied.

By Lemma~\ref{bcconstruction}, 
$b^2$ and $c^2$ have either a simple pole or a simple zero
at each branchpoint $\nu_1,\ldots,\nu_{2g}$. Since $\beta^2$ is
holomorphic on $\cnstar$ and vanishes at each of the points
$\nu_1,\ldots,\nu_{2g}$, $\beta^2b^2$ and $\beta^2c^2$ are holomorphic
functions on $\cnstar$. Furthermore, since by
Theorem~\ref{aconstruction}, $\beta^2a^2$ is the square of a
holomorphic function on $\cnstar$, the zero of $\beta^2$ at
each branchpoint is of odd order and $\beta^2$ can have no further
zeroes of odd order on $\cnstar$. Thus, $\beta^2b^2$ and
$\beta^2c^2$ have no zeroes of odd order, whence are squares of
holomorphic functions $\beta b$ and $\beta c$ on $\cnstar$.
\QED

\mytheorem{} {\em
Let $\FC$ be defined as in Section~\ref{FCconstruction}. Choose a
homotopy basis $a_1,\ldots,a_g,b_1,\ldots,b_g$ and the Abelian
differential $\Omega_1$ as in Section~\ref{FCconstruction}, such
that~\bref{qsuffcond} holds for some $q\in\cstar$.
Then there exists $0<r_0<1$ and 
a real $g$-parameter family $\FS_\FC$ of CMC-immersions,
such that each such immersion $\Psi\in\FS_\FC$ 
lies in the $r_0$-dressing orbit of the
cylinder, $q\in\Sym(\Psi)$, and the
hyperelliptic curve associated to $\Psi$ 
in Section~\ref{FCintro} is $\FC$. The family $\FS_\FC$ depends only
on $\FC$, not on the chosen solution of~\bref{qsuffcond}.
}

\Proof
By Theorem~\ref{aconstruction} and 
Proposition~\ref{bcconstruction}, to each hyperelliptic curve
$\FC$ and constant $q$, 
satisfying the conditions of Section~\ref{FCconstruction}, there
exists $0<r_0<1$, functions $f_+\in\FA_{r_0}^+$, $p=\lambda\inv
q-\lambda\qquer+f_+$, $\alpha=\cosh(p)$, $\beta=\sinh(p)$, $a^2$,
$b^2$, $c^2$, such that all assumptions of Theorem~\ref{sufftheorem}
are satisfied for $r^\prime=r=r_0$. 
{}From Theorem~\ref{aconstruction} we know that there exists at least a
$g$-parameter family of such data. Therefore, by the statement of
Theorem~\ref{sufftheorem}, there exists a real $g$-parameter family
$\FS_\FC$ of
CMC-immersions, such that each $\Psi\in\FS_\FC$ is
generated by dressing the cylinder with some
$h_+\in\Lambda_{r^\dprime}^+\LieSL(2,\Bcc)_\sigma$ with
$0<r^\dprime\leq r_0$ and $h_+Ah_+\inv=\tmatrix abc{-a}$. Since all
functions $a^2$ constructed in Theorem~\ref{aconstruction} are squares
of a meromorphic function $a$ on the same hyperelliptic curve $\FC$,
they all yield the surface $\FC$ under the construction in
Section~\ref{FCintro}.

The construction of $a^2$, $b^2$ and
$c^2$ in Sections~\ref{aconstruction} and~\ref{bcconstruction}, depends
only on $\FC$ not on the constant $q$. I.e.,
if there are two constants $q_1$, $q_2\in\cstar$ 
both satisfying~\bref{qsuffcond}, then for
both constants we get the same $0<r_0<1$ and $g$-parameter family
of matrices $h_+\in\Lambda_{r^\dprime}^+\LieSL(2,\Bcc)_\sigma$,
$r^\dprime\leq r_0$. Therefore, the family $\FS_\FC$ doesn't depend on
the chosen solution $q$ of~\bref{qsuffcond}.
\QED

\mycorollary{} {\em
Let $\FC$, $\Omega_1$, $q\in\cstar$ and $\Psi\in\FS_\FC$
be given as in Theorem~\ref{bcconstruction} and define $\omega$
by~\bref{pOmega}. For $\lambda_0\in S^1$ denote by $P_1(\lambda_0)$ and
$P_2(\lambda_0)$ the covering points of $\nu_0=\lambda_0^2$ on $\FC$.
Then the following are equivalent
\begin{description}
\item[1.] $q\in\Per(\Psi_{\lambda_0})$,
\item[2.] $\omega$ has a zero at $P_1(\lambda_0)$ or $P_2(\lambda_0)$ and
there exists a curve $\gamma$ on $\FCstar$ connecting $P_1(\lambda_0)$ and
$P_2(\lambda_0)$, such that $\int_\gamma\omega=2\pi im,\kern5mm m\in\Bii$.
\end{description}
}

\Proof
By Theorem~\ref{bcconstruction}, 
$q\in\Sym(\Psi)=\Sym(\Psi_\lambda)$ for all $\lambda\in S^1$, and the
curve defined for $\Psi$ in Section~\ref{FCintro} coincides with $\FC$.
Thus, the equivalence follows immediately from Theorem~\ref{necclosing}.
\QED

\newsection{} \label{Delaunay}
Theorem~\ref{bcconstruction} together with
Corollary~\ref{bcconstruction} allows us to classify periodic surfaces
in terms of algebro-geometric data: If $q\in\Per(\Psi_{\lambda_0})$,
$q\neq0$, then $\Per(\Psi_{\lambda_0})\neq\{0\}$, and, by
Theorem~\ref{intro11}, there exists a Riemann surface $M$ with
universal cover $\pi:\Bcc\rightarrow M$ 
and elementary group $\Gamma=\Per(\Psi_{\lambda_0})$
such that $\Psi_{\lambda_0}=\Phi\circ\pi$ defines a CMC-immersion
$\Phi:M\rightarrow\threespace$. In particular, if $\Per(\Psi_{\lambda_0})$ is
generated by only one translation, then $M$ is topologically a
cylinder. If there exist two linearly independent translations in
$\Per(\Psi_{\lambda_0})$, then $M$ is topologically a torus.
By part 3.\ in Theorem~\ref{intro11}, these are the only possible cases.
As a special case we get

\mytheorem{} {\em
Let $\FC$, $\Omega_1$ and $U_k$ be defined as in
Section~\ref{FCconstruction}. If there exists a real constant
$\phi\in[0,2\pi)$, such that 
\BEQ \label{Delaunaycond}
\Real(U_k)\sin(\phi)-\Imag(U_k)\cos(\phi)=0,\kern1cm k=1,\ldots,g,
\EEQ
then the associated family of each $\Psi$ constructed in
Theorem~\ref{bcconstruction} contains a Delaunay surface.

Conversely, if $\Psi(\Bcc)$ is in the associated family of a
Delaunay surface, and if $\FC$ and
$U_k$ are defined as in Sections~\ref{FCintro} and~\ref{pprime}, then
there exists $\phi\in[0,2\pi)$, such that~\bref{Delaunaycond} is satisfied.
}

\Proof
If $U_k$ is real, then for all real $q\in\Bcc$,~\bref{qsuffcond} is
satisfied with $m_k=0$, $k=1,\ldots,g$. Thus, $\Sym(\Psi)$ contains
the one-parameter subgroup $\Brr\subset\Bcc$. By Theorem~\ref{intro11} this
implies, that $\Psi$ is in the associated family of a Delaunay surface.

Conversely, if $\Psi(\Bcc)$ is a Delaunay surface, then,
by \cite[Lemma~2.15]{DoHa:2}, $\Sym(\Psi)$
contains a one-parameter subgroup generated by some $q_0\in\cstar$,
$|q_0|=1$. If we write $q_0=e^{i\phi}$ then~\bref{qsuffcond} holds for all 
$q$ of the form $q=re^{i\phi}$, $r\in\Brr$. Thus,
$\Imag(e^{i\phi}\overline{U_k})=0$ for all $k=1,\ldots,g$, which
gives~\bref{Delaunaycond}.
\QED

\newsection{} \label{sufftori}
Let us again consider the case of CMC-tori.

\mytheorem{} {\em
Let $\FC$, $\Omega_1$, $U_k$, $k=1,\ldots,g$, and
$q_1=q\in\cstar$ be given as in Theorem~\ref{bcconstruction}.
Assume that there
exists a second constant $q_2\in\cstar$, linearly independent of
$q_1$, such that also for $q_2$ Equation~\bref{qsuffcond} is
satisfied. Let $\lambda_0\in S^1$
and choose an arbitrary curve $\gamma$ in $\FCstar$ connecting the two points
$P_1(\lambda_0)$ and $P_2(\lambda_0)$, such that
$\nu(P_i(\lambda_0))=\lambda_0^2$. 
Define $X_1,Y_1,X_2,Y_2$, $\alpha_k,\beta_k$, $k=1,\ldots,g$,
$c_1,c_2$ as in Theorem~\ref{Bobenkotori}. If 1.\ and 2.\ of
Theorem~\ref{Bobenkotori} are satisfied, then there exists $0<r_0<1$
and a real $g-1$-parameter family $\FT_\FC$ of CMC-tori in the
$r_0$-dressing orbit of the cylinder, such that for each $\Psi\in\FT_\FC$, 
$q_1,q_2\in\Per(\Psi)$.
}

\Proof
Using the calculations in Section~\ref{Bobenkotori}, we see, that the
last $g$ columns of~\bref{intmatrix} are integer iff~\bref{qsuffcond}
holds for $q_1$ and for $q_2$. Therefore, by
Theorem~\ref{bcconstruction}, there exists a real $g$-parameter family
$\FS_\FC$ of CMC-immersions, such that
$q_1,q_2\in\Sym(\Psi)=\Sym(\Psi_\lambda)$ for all $\Psi\in\FS_\FC$ and all
$\lambda\in S^1$. Since all elements of the associated family of a
CMC-immersion have the same symmetry group, there is at least a
$g-1$-parameter family of associated families
$\{\Psi_\lambda;\lambda\in S^1\}$, such that $q_1$ and $q_2$ are in
$\Sym(\Psi_\lambda)$
for all $\lambda\in S^1$. From Corollary~\ref{bcconstruction} we get
that the first column of~\bref{intmatrix} is integer iff 
Condition~2. of Corollary~\ref{bcconstruction} is satisfied. By~\bref{dOI},
we have $\Omega_1(P_2)=-\Omega_1(P_1)$. Therefore,
$\Omega_1$ vanishes at $P_1$ iff it vanishes at $P_2$. Now we
define $\omega_1$ and $\omega_2$ for $q_1$ and $q_2$
by~\bref{pOmega}. If
$\Omega_1(P_1)=0$, then also, by~\bref{pOmega},
$\omega_1(P_1)=\omega_2(P_1)=0$. This shows, that for both, $q_1$ and $q_2$,
the conditions in Corollary~\ref{bcconstruction} are satisfied, i.e.,
$q_1,q_2\in\Per(\Psi_{\lambda_0})$. I.e., there exists a
$g-1$-parameter family $\FT_\FC=\{\Psi_{\lambda_0};\Psi=\Psi_1\in\FS_\FC\}$.
\QED

\separate
In~\cite{Ja:1}, existence of such families of
CMC-tori was proved for arbitrary genus $g>1$.

\myremark{}
1. In the case $g=1$, the differential $\Omega_1$ can be easily
computed explicitly: Since $\Omega_1$ has a pole of second order at
$P_\infty$, we know $\Omega_1=\frac{a\nu+b}{\mu}\diff\nu$. Writing
this in the coordinate $\lambda$ and comparing
with~\bref{dOsingbehaviour}, we see
$a=\frac12$. The function $b$ will be determined from 
the condition~\bref{dO1ak}. 
For this we represent the cycle $a_1$ by the straight line connecting
the points $\nu_1=r e^{i\phi}$ and $\nu_2=r\inv e^{i\phi}$. Then it is 
straightforward to compute
\BEQ
\int_{a_1}\Omega_1=e^{i\frac{\phi}{2}}\frac{2i}{\sqrt{r}}E(\sqrt{1-r^2})
+4ib\sqrt{r}e^{-i\frac{\phi}{2}}K(\sqrt{1-r^2}),
\EEQ
where $K$ and $E$ denote the complete elliptic differentials
of the first and second kind (see e.g.~\cite{GrRi:1}).
Now $\int_{a_1}\Omega_1=0$ yields
\BEQ
b=-\frac{e^{i\phi}}{2r}\frac{E(\sqrt{1-r^2})}{K(\sqrt{1-r^2})}.
\EEQ
If our choices would lead to a torus, then by
Theorem~\ref{Bobenkotori} we would
need to have $\Omega_1(\lambda_0)=0$ for some point $\lambda_0\in
S^1$. But this equivalent with $|b|=\frac12$ which, in view of our
formulas, is equivalent with
\BEQ
\frac{E(\sqrt{1-r^2})}{rK(\sqrt{1-r^2})}=1.
\EEQ
But the left hand side is always $>1$ for $0<r<1$, producing a contradiction.

2. If we only want to produce Delaunay surfaces from hyperelliptic
curves of genus $1$, then we start as above and consider
$\Omega_1$. The coefficients $a$ and $b$ are determined as
above. However, we do not need to require $\Omega_1(P_1(\lambda_0))=0$
for some $\lambda_0\in S^1$ as above, but
only~\bref{Delaunaycond}. Since $g=1$, this is only one equation for
$U_1$. Obviously, this equation is solvable. Thus in the case $g=1$
our construction produces a family of Delaunay surfaces,
indexed by $0<r<1$ and the scale factor of $a^2$ in our construction.
For a thorough investigation of the case $g=1$ see~\cite[Chapter~8]{PiSt:1}.

3. It is well known, that there are only two one-parameter families of Delaunay
surfaces, the unduloids and the nodoids~\cite{Ee:1}. Thus, the
discussion of Delaunay surfaces above also shows, that the
$g$-parameter families constructed above can contain many congruent
surfaces. It would certainly be interesting to further investigate
this question.

\newsection{} \label{interpretation}
In Section~\ref{tori} we have reproduced the algebro-geometric 
description of CMC-tori~\cite{Bo:1,PiSt:1} using the $r$-dressing method.
It is important to note, that our starting point differs from the one
in~\cite{PiSt:1} and~\cite{Bo:1}. 

Pinkall and Sterling started with the observation, that the metric of
a CMC-torus is a finite type solution of the $\sinh$-Gordon equation.
Thus, as follows from the Krichever construction of finite type
solutions to integrable systems, there is a hyperelliptic curve
associated to each CMC-torus. However, the analytic construction of
this surface from a finite type solution is not unique. Actually, for
each trivial flow of a given finite type solution, there is a, possibly
singular, hyperelliptic curve. The ambiguity of the definition of the
hyperelliptic curve introduces an additional level of complication in
the analytic description. While Pinkall and Sterling restricted their
attention to nonsingular hyperelliptic curves,
Bobenko~\cite[Appendix]{Bo:1} gave an argument why CMC-tori belong
generically to nonsingular curves.

Our approach is more geometric in the sense, that we start directly from
the periodicity conditions of the extended frame of a CMC-immersion.
In Section~\ref{periodicsurfaces}, we introduced a hyperelliptic curve
$\FC$, which is directly derived from the periodicity conditions on
the extended frame of a CMC-immersion with periodic metric. By
Theorem~\ref{finitetype}, if the metric of a CMC-immersion in the
$r$-drssing orbit of the standard cylinder is periodic, then it is of
finite type, and our hyperelliptic curve $\FC$ coincides with one
of those used in the Krichever construction. However, by definition, the
hyperelliptic curve chosen in this paper 
is always nonsingular. We thereby circumvent the
problem of dealing with singular curves.

Of course, there still remains the question, if the construction in
Sections~\ref{FCconstruction}--\ref{bcconstruction} gives all CMC-tori.

Let $\Psi:\Bcc\rightarrow\threespace$ be a
CMC-torus. Then, by \cite[Corollary~5.3]{DoWu:1}, 
$\Psi$ is in the $r$-dressing orbit
of the standard cylinder for some $0<r<1$.
Hence, by Theorem~\ref{nectheorem}, there are rational functions
$a^2$, $b^2$ and $c^2$, satisfying certain conditions,
such that $\Psi$ is obtained by $r$-dressing
with $h_+\in\Lambda_r\LieSL(2,\Bcc)_\sigma$, where
$h_+Ah_+\inv=\tmatrix abc{-a}$.
Moreover, by Theorem~\ref{Bobenkotori}, there is a nonsingular
hyperelliptic curve $\FC$ such that $a^2$ is the square of a
meromorphic function on $\FC$. 

Conversely, if we start from a hyperelliptic curve $\FC$ as in
Section~\ref{FCconstruction}, then the construction in
Section~\ref{aconstruction} and Section~\ref{bcconstruction} 
gives only functions $a^2$, $b^2$ and $c^2$ with simple poles.
It is easy to check, that the construction gives in fact all rational
functions $a^2,b^2,c^2$ with only simple poles, which have the
properties stated in Theorem~\ref{nectheorem}.

The question, if this is enough to get all CMC-tori, leads to a
statement analogous to the one made in the appendix of~\cite{Bo:1}:

\myconjecture{}
Let $\Psi:\Bcc\rightarrow\threespace$ be a CMC-immersion, which is
obtained by dressing the cylinder with
$h_+\in\Lambda_r^+\LieSL(2,\Bcc)_\sigma$ for some $0<r<1$. Define
$h_+Ah_+\inv=\tmatrix abc{-a}$. If $\Psi(\Bcc)$ is a CMC-torus, then
the function $a^2(\lambda)$ is rational and has generically only
poles of first order.

\separate
We would think that this is likely to be true. In this case,
generically, to each hyperelliptic curve $\FC$
there would exist either no CMC-torus or, up to finitely many choices,
precisely the $g-1$-parameter
family of CMC-tori $\FT_\FC$ constructed in Theorem~\ref{sufftori}.

We would expect that the understanding of singular tori ---if they
exist--- will provide the technical tools to answer these questions.

\end{document}